\documentstyle[amsfonts, amssymb, latexsym, epic, eepic, 12pt]{amsart}

\newtheorem{theorem}{Theorem}[section]
\newtheorem{lemma}[theorem]{Lemma}
\newtheorem{proposition}[theorem]{Proposition}
\newtheorem{corollary}[theorem]{Corollary}
\newtheorem{definition}[theorem]{Definition}

\def\C{{\mbox{\rm\kern.24em
\vrule width.03em height1.43ex depth-.052ex \kern-.26em C}}}
\def\QSet{\mbox{\rm\kern.24em
\vrule width.03em height1.48ex depth-.051ex \kern-.26em Q}}
\def\Z{{\bf Z}}
\def\R{{\mbox{\rm I\kern-.22em R}}}
\def\P{{\bf P}}
\def\Q{{\bf Q}}
\def\T{{\bf T}}
\def\D{{\bf D}}

\def\size{{\rm size}}

\def\energy{{\rm energy}}

\def\M{{\rm M}}

\def\111{\gamma}

\def\be#1{\begin{equation}\label{#1}}
\def\bas{\begin{align*}}
\def\eas{\end{align*}}
\def\bi{\begin{itemize}}
\def\ei{\end{itemize}}
\newenvironment{proof}{\noindent {\bf Proof} }{\endprf\par}
\def \endprf{\hfill  {\vrule height6pt width6pt depth0pt}\medskip}
\def\emph#1{{\it #1}}

\title{On the bi-carleson operator I. The Walsh case}

\author{Camil Muscalu}
\address{Department of Mathematics, UCLA, Los Angeles CA 90095-1555}
\email{camil@@math.ucla.edu}

\author{Terence Tao}
\address{Department of Mathematics, UCLA, Los Angeles CA 90095-1555 }
\email{tao@@math.ucla.edu}

\author{Christoph Thiele}
\address{Department of Mathematics, UCLA, Los Angeles CA 90095-1555}
\email{thiele@@math.ucla.edu}

\include{psfig}
\begin{document}

\begin{abstract}  
We prove $L^p$ estimates 
for the Walsh model of the maximal bi-Carleson operator defined below.  
The corresponding estimates for the Fourier model will be obtained in 
the sequel \cite{mtt:fourierbicarleson} of this paper.
\end{abstract}

\maketitle

\section{introduction}

In this paper we initiate a study of the bi-Carleson maximal operator which, as we will see, is a hybrid of the Carleson maximal operator and
the bilinear Hilbert transform.

Let us first recall that the maximal Carleson operator is the 
sub-linear operator defined by

\begin{equation}
C(f)(x):=\sup_{N}\left|\int_{\xi < N}
\widehat{f}(\xi) e^{2\pi ix\xi}\, d\xi\right|,
\end{equation}
where $f$ is a Schwartz function on the real line $\R$ and the Fourier 
transform is defined by 

$$ \hat f(\xi) := \int_\R  f(x)e^{-2\pi i x \xi}\ dx.$$
The following statement of Carleson and Hunt 
\cite {carleson}, \cite{hunt} is a classical theorem 
in Fourier analysis:

\begin{theorem}
\cite{carleson}, \cite{hunt}
The operator $C$ maps $L^p\rightarrow L^p$, for every $1<p<\infty$.
\end{theorem}
This result, in the particular weak type $(2,2)$ special case,
 was the main ingredient in the proof of Carleson's fameous
theorem which states that the Fourier series of a function in
$L^2[0,1]$ converges pointwise almost everywhere.

The bilinear Hilbert transform is an operator which can be essentially written  as
\begin{equation}
B(f_1,f_2)(x):=\int_{\xi_1<\xi_2}
\widehat{f}_1(\xi_1)\widehat{f}_2(\xi_2) e^{2\pi ix(\xi_1+\xi_2)}\,d\xi_1 d\xi_2,
\end{equation}
where $f_1, f_2$ are test functions on $\R$.

From the work of Lacey and Thiele \cite{laceyt1}, \cite{laceyt2} we have the following $L^p$ estimates on $B$:

\begin{theorem}\label{bht}\cite{laceyt1}, \cite{laceyt2}
$B$ maps $L^p \times L^q \to L^r$ whenever $1 < p, q \leq \infty$, $1/p + 1/q = 1/r$, and $2/3 < r < \infty$.
\end{theorem}

In this paper and the sequel \cite{mtt:fourierbicarleson} we study the $L^p$ boundedness properties of a 
maximal sub-bilinear variant $T$ of the bilinear Hilbert transform,
 defined by \footnote{The operator $T_-$ defined by the same formula but with phase  
 "`$\xi_1-\xi_2$"' instead of "`$\xi_1+\xi_2$"', has recently been studied in  \cite{mtt:counterexample}: it does not satisfy any $L^p$ estimates!}
\begin{equation}\label{oper}
T(f_1,f_2)(x):=\sup_{N}\left|\int_{\xi_1<\xi_2< N}
\widehat{f}_1(\xi_1)\widehat{f}_2(\xi_2)
e^{2\pi ix(\xi_1+\xi_2)}\,d\xi_1 d\xi_2\right|.
\end{equation}
As we already discussed in \cite{mtt:walshbiest}, \cite{mtt:fourierbiest},
operators like $T$ are closely related to the 
WKB expansions of generalized eigenfunctions of 
one-dimensional Schr\"{o}dinger operators, 
following the work of Christ and Kiselev \cite{ck}.

The study of $T$ clearly reduces to the study of its linearized version
defined by

\begin{equation}\label{bicarleson-def}
T(f_1,f_2)(x):=\int_{\xi_1<\xi_2< N(x)}
\widehat{f}_1(\xi_1)\widehat{f}_2(\xi_2)
e^{2\pi ix(\xi_1+\xi_2)}\,d\xi_1 d\xi_2,
\end{equation}
where ``$x\rightarrow N(x)$'' is an arbitrary function, 
in the same way  as the study of $C$ reduces to the study of 

\begin{equation}\label{carleson-def}
C(f)(x):=\int_{\xi < N(x)}
\widehat{f}(\xi) e^{2\pi ix\xi}\, d\xi.
\end{equation}
From the identity
$$
f_1(x) f_2(x)  = \int
\widehat{f}_1(\xi_1)\widehat{f}_2(\xi_2)
e^{2\pi ix(\xi_1+\xi_2)}\,d\xi_1 d\xi_2 
$$
we see that $T$ has the same homogeneity as the pointwise 
product operator, and hence we expect estimates of H\"older type, i.e. that $T$ maps $L^{p_1} \times L^{p_2}$  to $L^{p'_3}$ when $1/p'_3= 1/p_1 + 1/p_2$.

We should emphasize that the relationship between the bi-Carleson
operator $T$ and the bilinear
Hilbert transform $B$, is analogous to the relationship between the Carleson
operator $C$ and the classical Hilbert transform.\footnote{This should also be compared with
M.Lacey's maximal operator \cite{lacey}, which in this sub-bilinear context is the analogue of
the Hardy Littlewood maximal function.}

It is well-known that the Carleson operator $C$ and the bilinear Hilbert transform $B$ have slightly simpler 
Walsh model versions $C_{walsh}$ and $B_{walsh}$ defined using the Walsh 
transform instead of the Fourier transform, which we will present next.
We should also recall here that here are several papers in the literature using these Walsh notations
(see e.g. \cite{billard}, \cite{thiele}, \cite{mtt:walshbiest}, etc.).

\begin{definition}\label{walsh-def}
For $l\geq 0$ we define the $l$-th \emph{Walsh function}
$w_l$ by the following recursive formulas

\begin{eqnarray*}
w_0         & := & \chi_{[0,1)}         \\
w_{2l}      & := & w_l(2x)+w_l(2x-1)    \\
w_{2l+1}    & := & w_l(2x)- w_l(2x-1).  
\end{eqnarray*}
\end{definition}

\begin{definition}\label{tile-def}
A \emph{tile} $P$ is a half open rectangle $I_P\times\omega_P$ of area one, such that
$I_P$ and $\omega_P$ are dyadic intervals. If 
$P=[2^{-k}n, 2^{-k}(n+1))\times [2^k l, 2^k(l+1))$ 
is such a tile, we
define the corresponding \emph{Walsh wave packet} $\phi_P$ by
\[\phi_P(x):= 2^{k/2}w_l(2^k x-n).\]
\end{definition}

For each tile $P$, note that $\phi_P$ is supported on $I_P$ and has an $L^2$ norm equal to 1.  Also, observe that $\phi_P$ and $\phi_{P'}$ are orthogonal whenever $P$ and $P'$ are disjoint.

\begin{definition}\label{bitile-def}
A \emph{bitile} $P$ is a half open rectangle $I_P\times\omega_P$ of area two,
 such that
$I_P$ and $\omega_P$ are dyadic intervals . 
For any bitile 
\[ P = [2^{-k}n, 2^{-k}(n+1))\times [2^{k+1} l, 2^{k+1}(l+1))\]
we define the sub-tiles $P_1, P_2\subset P$ by
\begin{eqnarray*}
P_1 & := & [2^{-k}n, 2^{-k}(n+1))\times [2^k 2l, 2^k(2l+1)) \\
P_2 & := & [2^{-k}n, 2^{-k}(n+1))\times [2^k(2l+1), 2^k(2l+2)). \\
\end{eqnarray*}
\end{definition}

\begin{definition}\label{quartile-def}
A \emph{quartile} $P$ is a half 
open rectangle $I_P\times\omega_P$ of area four, 
such that
$I_P$ and $\omega_P$ are dyadic intervals . 
For any quartile 
\[ P = [2^{-k}n, 2^{-k}(n+1))\times [2^{k+2} l, 2^{k+2}(l+1))\]
we define the sub-tiles $P_1, P_2, P_3 \subset P$ by
\begin{eqnarray*}
P_1 & := & [2^{-k}n, 2^{-k}(n+1))\times [2^k 4l, 2^k(4l+1)) \\
P_2 & := & [2^{-k}n, 2^{-k}(n+1))\times [2^k(4l+1), 2^k(4l+2)) \\
P_3 & := & [2^{-k}n, 2^{-k}(n+1))\times [2^k(4l+2), 2^k(4l+3)).
\end{eqnarray*}
and the sub-bitile $P_{12}$ by $P_{12}=P_1\cup P_2$.
\end{definition}
Note that every quartile is a disjoint union of two bitiles and also
a disjoint union of four tiles.

\begin{definition}\label{wcarleson-def}(\cite{billard},\cite{sjolin})
If $\P$ is a finite collection of bitiles, the Walsh Carleson operator
is defined by the formula
$$C_{walsh,\P}(f):=\sum_{P\in\P}
\langle f,\phi_{P_1}\rangle
\phi_{P_1}\chi_{\{x:N(x)\in\omega_{P_2}\}}.$$
\end{definition}
Similarly, one has

\begin{definition}\label{bht-def}(\cite{thiele})  If $\P$ is a finite collection of quartiles, the Walsh bilinear Hilbert transform $B_{walsh,\P}$ is defined by the formula
$$ B_{walsh,\P}(f_1,f_2) := \sum_{P \in \P} \frac{1}{|I_P|^{1/2}}
\langle f_1, \phi_{P_1} \rangle \langle f_2, \phi_{P_2} \rangle \phi_{P_3}.$$
\end{definition}

From the point of view of the time-frequency phase plane, the Fourier and Walsh models are very similar\footnote{Indeed, the Walsh model can be viewed as the analogue of the Fourier model with the underlying group $\R$ being replaced by $(\Z_2)^{\Z}$.}, however the Walsh model, being dyadic, has several convenient features (such as the ability to localize perfectly in both time and frequency simultaneously) which allow for a clearer and less technical treatment than the Fourier case.  
The following theorems are well known:

\begin{theorem}\label{carleson-walsh}\cite{billard},\cite{sjolin}
Let $\P$ be an arbitrary set of bitiles. Then,
$C_{walsh,\P}$ maps $L^p\rightarrow L^p$ for every $1<p<\infty$.
The bounds are uniform in $\P$
\end{theorem}

\begin{theorem}\label{bht-walsh}\cite{thiele'} 
Let $\P$ be an arbitrary set of quartiles.
$B_{walsh,\P}$ maps $L^p \times L^q \to L^r$ whenever $1 < p, q \leq \infty$, $1/p + 1/q = 1/r$, and $2/3 < r < \infty$.  The bounds are uniform in $\P$.
\end{theorem}

Just as the Carleson operator $C$ and the bilinear Hilbert transform $B$ 
have Walsh models $C_{walsh,\P}$ and $B_{walsh,\P}$, the operator $T$ also 
has a Walsh model $T_{walsh,\P,\Q}$. 

\begin{definition}\label{twalsh-def}  If $\P$, $\Q$ are two finite collections of quartiles, we define the operator $T_{walsh,\P,\Q}$ by
$$T_{walsh,\P,\Q}:=T'_{walsh,\P,\Q}+T''_{walsh,\P,\Q}$$
where
\bas
T'_{walsh,\P,\Q}(f_1,f_2)&:=
\sum_{P\in\P}
\langle B_{walsh,P_1,\Q}(f_1,f_2),\phi_{P_1}\rangle
\phi_{P_1}\chi_{\{x:N(x)\in\omega_{P_2}\}}  \\
T''_{walsh,\P,\Q}(f_1,f_2)&:=
\sum_{P\in\P}
\langle f_1,\phi_{P_1}\rangle
\phi_{P_1}C_{walsh,P_2,\Q}(f_2)\chi_{\{x:N(x)\in\omega_{P_2}\}}
\end{align*}
where for every tile $P$, $B_{walsh,P,\Q}$ and $C_{walsh,P,\Q}$ are defined by
\bas
B_{walsh,P,\Q}(f_1,f_2)&:=
\sum_{Q\in\Q;\, \omega_{Q_3}\subseteq \omega_P}
\frac{1}{|I_Q|^{1/2}}
\langle f_1,\phi_{Q_1}\rangle
\langle f_2,\phi_{Q_2}\rangle
\phi_{Q_3} \\
C_{walsh,P, \Q}(f_2)&:=
\sum_{Q\in\Q;\, \omega_{Q_1}\subseteq \omega_P}
\langle f_2,\phi_{Q_1}\rangle
\phi_{Q_1}\chi_{\{x:N(x)\in\omega_{Q_2}\}}.
\end{align*}
\end{definition}

In the next section we will explain why the operator $T_{walsh,\P,\Q}$ is the natural Walsh analogue of the Fourier operator $T$.  

The purpose of this paper is to obtain a large set of $L^p$ estimates for the Walsh model $T_{walsh,\P,\Q}$ of $T$.  The operator $T$ itself is 
 more technical to handle, and the treatment will be deferred to the sequel
 \cite{mtt:fourierbicarleson} of this paper.  
Our main theorem is the following:

\begin{theorem}\label{main}
Let $\P$, $\Q$ be two arbitrary finite collections of quartiles. Then,
$T_{walsh,\P,\Q}$ maps
\begin{equation}\label{star}
T_{walsh,\P,\Q}: L^{p_1}\times L^{p_2}\rightarrow L^{p'_3}
\end{equation}
as long as $1<p_1, p_2\leq\infty$, 
$\frac{1}{p_1}+\frac{1}{p_2}=\frac{1}{p'_3}$ and $2/3<p'_3<\infty$,
with uniform bounds in $\P$ and $\Q$.
\end{theorem}

We should point out that both Theorem \ref{carleson-walsh} and
Theorem \ref{bht-walsh} are particular cases of our main Theorem \ref{main}.

Also, as a consequence of our framework, we shall be able to give a self-contained proof of theorem
\ref{carleson-walsh} in Section \ref{carleson-walsh-sec}, in the same spirit with the proof of theorem
\ref{bht-walsh} in \cite{mtt:walshbiest}.

While the current article is mostly selfcontained, we adopt the same
attitude as in \cite{mtt:walshbiest}, \cite{mtt:fourierbiest} and will mark
as ``standard'' any arguments that are well understood by now in the framework
of multilinear singular integrals as in
\cite{laceyt1}, \cite{laceyt2}, \cite{cct}, 
\cite{mtt:walshbiest}, \cite{mtt:fourierbiest},
\cite{thiele'}, \cite{thiele}, etc.

The first author was partially supported by the NSF grant DMS 0100796.
The second author is a Clay Prize Fellow and is supported by grants from
NSF and Packard Foundation. The third author was partially supported
by a Sloan Fellowship and by NSF grants DMS 9985572 and DMS 9970469.

\section{the walsh model}

We explain in this section the analogy between $T_+$ and $T_{walsh,\P,\Q}$.
Here $T_+$ is defined in the same way as $T$,
 but with respect to an integration
over $0<\xi_1<\xi_2$ instead. This is a harmless modification (essentially
one replaces $f_1$ with the Hilbert transform of $f_1$).
Formally, 
one can write

\begin{equation}\label{unauna}
T_+(f_1,f_2)=
\sum_{\omega}
\int_{\xi_1<\xi_2<N(x);\,\xi_1,\xi_2\in\omega_{l};\,
N(x)\in\omega_{r}}
\widehat{f_1}(\xi_1)\widehat{f_2}(\xi_2)
e^{2\pi ix(\xi_1+\xi_2)}\,d\xi_1 d\xi_2 +
\end{equation}

\[\sum_{\omega}\int_{\xi_1<\xi_2<N(x);\,\xi_1\in\omega_{l};\,
\xi_2, N(x)\in\omega_{r}}
\widehat{f_1}(\xi_1)\widehat{f_2}(\xi_2)
e^{2\pi ix(\xi_1+\xi_2)}\,d\xi_1 d\xi_2.\]

This is a consequence of the fact that for almost every $3$-tuple
$\xi_1<\xi_2<N(x)$ there is a unique smallest dyadic interval
$\omega$ which contains $\xi_1,\xi_2$ and $N(x)$ and either
$\xi_1,\xi_2\in\omega_{l}$ and $N(x)\in\omega_{r}$, or
$\xi_1\in\omega_{l}$ and $\xi_2,N(x)\in\omega_{r}$, where
$\omega_{l}(\omega_{r})$ is the left (right) half of $\omega$.

Then (\ref{unauna}) is equal to

\begin{equation}\label{lakers}
\sum_{\omega}\left(\int_{\xi_1<\xi_2;\,\xi_1,\xi_2\in\omega_{l}}
\widehat{f_1}(\xi_1)\widehat{f_2}(\xi_2)
e^{2\pi ix(\xi_1+\xi_2)}\,d\xi_1 d\xi_2 \right)
\chi_{\{ x:N(x)\in\omega_r\} }+
\end{equation}

\begin{equation}\label{chicago}
\sum_{\omega}\left(\int_{\xi_1\in\omega_{l}}\widehat{f}_1(\xi_1)
e^{2\pi ix\xi_1}\,d\xi_1\right)
\left(\int_{\xi_2<N(x);\,\xi_2,N(x)\in\omega_{r}}
\widehat{f_2}(\xi_2)
e^{2\pi ix(\xi_2)}\,d\xi_2 \right).
\end{equation}
In the end, one just has to observe that the left hand side
term in (\ref{lakers}) is the bilinear Hilbert transform
of $f_1$ and $f_2$ whose Fourier transform is restricted to $\omega_l$, while
the right hand side term in (\ref{chicago}) is a Carleson operator restricted to
$\omega_r$.

\section{interpolation}\label{interp-sec}

In this section we review the interpolation theory from \cite{cct} which 
allows us to reduce multi-linear $L^p$ estimates such as those in 
Theorem \ref{main} to certain ``restricted  type'' estimates.

Throughout the paper, we use $A\lesssim B$ to denote the statement that
$A\leq CB$ for some large constant $C$, and $A\ll B$ to denote the
statement that $A\leq C^{-1}B$ for some large constant $C$. 
Also, the function $N(x)$ is fixed throughtout the paper and
our constants $C$ will always be independent of $\P$, $\Q$ and $N(x)$.

To prove the $L^p$ estimates on $T_{walsh,\P,\Q}$ it is convenient to use 
duality and introduce the trilinear form $\Lambda_{walsh,\P,\Q}$ associated to $T_{walsh,\P,\Q}$ via the formula
$$
\Lambda_{walsh,\P,\Q}(f_1, f_2,f_3) := \int_{\R} 
T_{walsh,\P,\Q}(f_1, f_2)(x)f_3(x) dx.
$$
Similarly define $\Lambda'_{walsh,\P,\Q}$ and $\Lambda''_{walsh,\P,\Q}$.  The statement that $T_{walsh,\P,\Q}$ is bounded from 
$L^{p_1} \times L^{p_2}$ to $L^{p'_3}$ is then equivalent to 
$\Lambda_{walsh,\P,\Q}$ being bounded on 
$L^{p_1} \times L^{p_2} \times L^{p_3}$ if 
$1 \le  p_3' < \infty$.  For $p_3'<1$ this simple duality relationship breaks down, however the interpolation arguments in \cite{cct} will allow us to reduce (\ref{star}) to certain ``restricted
type'' estimates on $\Lambda_{walsh,\P,\Q}$. As in \cite{cct} we find more
convenient to work with the quantities $\alpha_i=1/p_i$, $i=1,2,3$,
where $p_i$ stands for the exponent of $L^{p_i}$.

\begin{definition}
A tuple $\alpha=(\alpha_1,\alpha_2,\alpha_3)$ is called admissible, if

\[-\infty <\alpha_i <1\]
for all $1\leq i\leq 3$,

\[\sum_{i=1}^3\alpha_i=1\]
and there is at most one index $j$ such that $\alpha_j<0$. We call
an index $i$ good if $\alpha_i\geq 0$, and we call it bad if
$\alpha_i<0$. A good tuple is an admissible tuple without bad index, a
bad tuple is an admissible tuple with a bad index.
\end{definition}

\begin{definition}
Let $E$, $E'$ be sets of finite measure. We say that $E'$ is a major 
subset of $E$ if $E'\subseteq E$ and $|E'|\geq\frac{1}{2}|E|$.
\end{definition} 

\begin{definition}
If $E$ is a set of finite measure, we denote by $X(E)$ the space of
all functions $f$ supported on $E$ and such that $\|f\|_{\infty}\leq
1$.
\end{definition}

\begin{definition}
Let $\alpha=(\alpha_1,\alpha_2,\alpha_3)$ be an admissible tuple.
If $\alpha$ is a bad tuple, let $j$ be its bad index; otherwise let $j$ be
arbitrary. We say
that a $3$-linear form $\Lambda$ is of restricted type $\alpha$
if for every sequence $E_1, E_2, E_3$ of subsets of $\R$ with
finite measure, there exists major subsets $E'_i$ of $E_i$ such that

\[|\Lambda(f_1, f_2, f_3)|\lesssim |E|^{\alpha}\]
for all functions $f_i\in X(E'_i)$, $i=1,2,3$, where we adopt the
convention $E'_i =E_i$ for indices $i\neq j$, and $|E|^{\alpha}$
is a shorthand for

\[|E|^{\alpha}=|E_1|^{\alpha_1}|E_2|^{\alpha_2}
|E_3|^{\alpha_3}.\]

\end{definition}

Let us consider now the $2$-dimensional affine hyperspace
\[
S:=\{(\alpha_1,\alpha_2,\alpha_3)\in\R^3\,
|\,\alpha_1 + \alpha_2 + \alpha_3=1\}.
\]

\begin{figure}[htbp]\centering
\psfig{figure=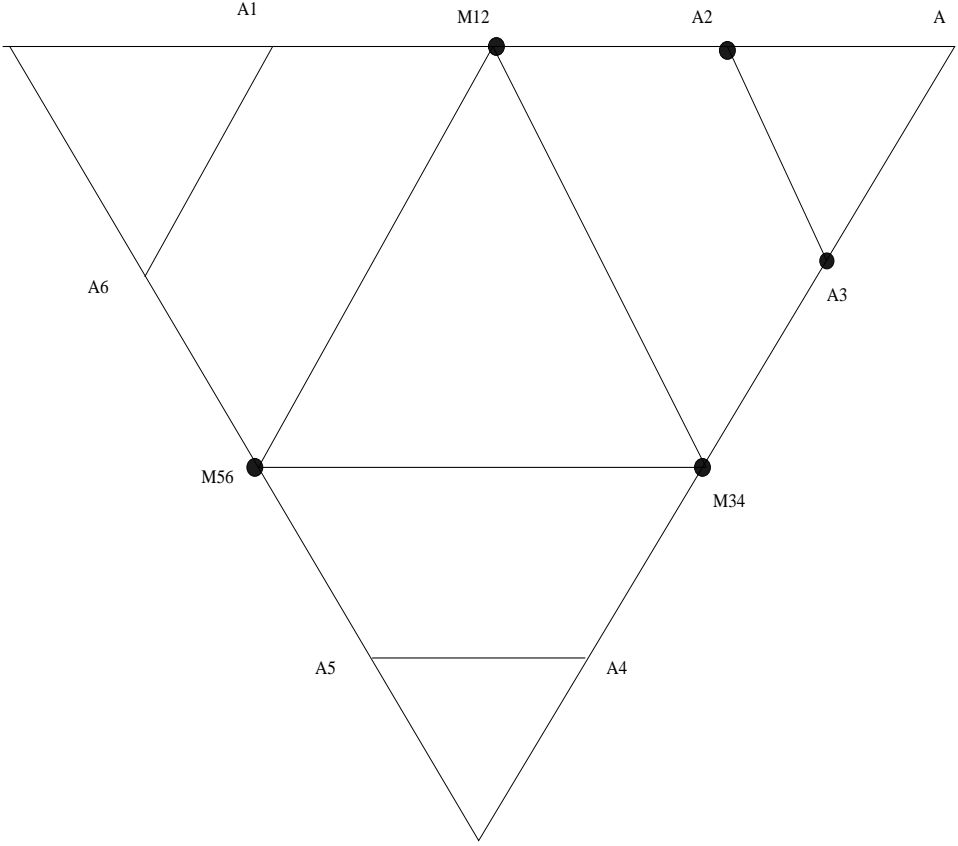,height=4.4in,width=5in}
\caption{The polytope $[A_1...A_{6}]$}
\label{fig}
\end{figure}

The points $A_1,...,A_6$ belong to $S$ and have the following coordinates:

\[
\begin{array}{llll}
A_1:(-\frac 12,1,\frac 12)  &  A_2:(\frac 12,1,-\frac 12)  &  
A_3:(1,\frac 12,-\frac 12)   \\  
\ &\ &\ & \ \\
A_4:(1,-\frac 12,\frac 12)  &  A_5:(\frac 12,-\frac 12,1) & 
A_6:(-\frac 12,\frac 12,1).   \\  
\end{array}
\]
The points $M_{12}, M_{34}, M_{56}$ are midpoints of their
corresponding segments and have the coordinates
$M_{12}:(0,1,0)$, $M_{34}:(1,0,0)$, $M_{56}:(0,0,1)$.

The following ``restricted type'' results will be proved directly.

\begin{theorem}\label{teorema1}
For every vertex $A_i$, $i = 1, \ldots, 6$ there exist
admissible tuples $\alpha$ arbitrarily close to $A_i$ such that the 
form $\Lambda'_{walsh,\P,\Q}$ is of restricted
type $\alpha$ uniformly in $\P$, $\Q$.
\end{theorem}

\begin{theorem}\label{teorema2}
For the vertices $M_{56}, M_{12}, A_2$ there exist
admissible tuples $\alpha$ arbitrarily close to them, such that the 
form $\Lambda''_{walsh,\P,\Q}$ is of restricted
type $\alpha$ uniformly in $\P$, $\Q$.
\end{theorem}

By theorem \ref{teorema1} and the interpolation theory in 
\cite{cct} Section 3, we obtain

\begin{corollary}\label{restricted1}
Let $\alpha$ be an admissible tuple inside the hexagon $[A_1,...,A_6]$.
Then $\Lambda'_{walsh,\P,\Q}$ is of restricted type
$\alpha$.
\end{corollary}
We also have

\begin{corollary}\label{restricted2}
Let $\alpha$ be an arbitrary tuple inside the pentagon
$[M_{56} M_{34} A_3 A_2 M_{12}]$.
Then $\Lambda''_{walsh,\P,\Q}$ is of restricted type
$\alpha$.
\end{corollary}

\begin{proof}
By using the interpolation theory in \cite{cct} and theorem \ref{teorema2},
it is enough to show that there exist admissible tuples $\alpha$
arbitrarily close to $M_{34}$ and $A_3$, such that the form
$\Lambda''_{walsh,\P,\Q}$ is of restricted type $\alpha$.

But one can see this immediately from the identity:

$$\Lambda''_{walsh,\P,\Q}(f_1,f_2)=C_{walsh,\P}(f_1)\cdot C_{walsh,\Q}(f_2)-
\Lambda''_{walsh,\Q,\P}(f_2,f_1)$$
which can be obtained from definition \ref{twalsh-def}, after changing
the order of summation over $\P$ and $\Q$.

Indeed, one then just has to apply theorem \ref{teorema2}
to handle the form $\Lambda''_{walsh,\Q,\P}(f_2,f_1)$ and theorem
\ref{carleson-walsh} to handle the product
$C_{walsh,\P}(f_1)\cdot C_{walsh,\Q}(f_2)$.
\end{proof}

Intersecting these two corollaries we obtain the analogous result for 
$\Lambda$ and the pentagon $[M_{56} M_{34} A_3 A_2 M_{12}]$.

Since one observes that $p_1, p_2, p'_3$ satisfy the hypothesis
of theorem \ref{main} if and only if 
$(1/p_1,1/p_2,1/p_3)\in [M_{56} M_{34} A_3 A_2 M_{12}]$,
it only remains to convert these restricted type estimates into strong type estimates.  
To do this, one just has to apply (exactly as in \cite{cct}) the multilinear Marcinkiewicz interpolation theorem \cite{janson} in the case of good tuples and the interpolation lemma 3.11 in \cite{cct} in the case of bad tuples.

This ends the proof of theorem \ref{main}. Hence, it remains to prove
theorem \ref{teorema1} and theorem \ref{teorema2}.

\section{trees}

In order to prove the desired estimates for the forms $\Lambda'_{walsh,\P,\Q}$
and $\Lambda''_{walsh,\P,\Q}$ one needs to organize our collections of
quartiles $\P$, $\Q$ into trees as in \cite{fefferman}, 
\cite{laceyt1}, \cite{laceyt2}, \cite{cct}.

\begin{definition}
Let $P$ and $P'$ be tiles. We write $P'<P$ if $I_{P'}\subsetneq I_P$
and $\omega_P\subseteq \omega_{P'}$, and $P'\leq P$ if $P'<P$ or
$P'=P$. 
\end{definition}

Note that $<$ forms a partial order on the set of tiles.  
Sometime, we will use the same partial order for bi-tiles (see 
Definition \ref{size'energy'-def}).
The transitivity of this order shall be crucial, especially in Lemma \ref{biest-trick}.

\begin{definition}
For every $j=1,2,3$ and $P_T\in\P$ define a $j$-tree with top
$P_T$ to be a collection of quartiles $T\subseteq\P$ such that

\[P_j\leq P_{T,j}\]
for all $P\in T$. 
We also say that $T$ is a tree if it is
a $j$-tree for some $j=1,2,3$.
\end{definition}

Notice that $T$ does not necessarily have to contain its top $P_T$.

The following geometric lemma is standard and easy to prove (see 
\cite{laceyt1}, \cite{laceyt2}, \cite{cct}).

\begin{lemma}\label{lacunar}
Let $P, P'$ be quartiles, and $i,j = 1,2,3$ be such that $i \neq j$. If
$P'_i\leq P_i$ then $P'_j\cap P_j=\emptyset$.

In particular, if $T$ is an $i$-tree, then the tiles $\{ P_j: P \in T\}$ are pairwise disjoint.
\end{lemma}

Even more is true: if $T$ is an $i$- tree, then the elements $P$ in
$T$ are parameterized by $I_P$ and the functions $\phi_{P_j}$ behave
like Haar functions in the sense that Calderon- Zygmund theory applies.
Thus an $i$-tree $T$ may be called lacunary in the indices $j$ other 
than $i$. 

\section{Tile norms, part 1}

In this section we start the study of the form $\Lambda'_{walsh,\P,\Q}$.
To do this, we shall frequently estimate expressions of the form
\be{trilinear}
| \sum_{P \in \P} \frac{1}{|I_P|^{1/2}} a^{(1)}_{P_1} 
a^{(2)}_{P_2} a^{(3)}_{P_3}|
\end{equation}
or of the form

\begin{equation}\label{bilinear}
| \sum_{P\in\P}a_{P_1} b_{P_2}|
\end{equation}
where $\P$ is a collection of quartiles, $a^{(j)}_{P_j}$ are complex 
numbers of the type

\begin{equation}
 a^{(j)}_{P_j} = \langle F_P, \phi_{P_j} \rangle
\end{equation}
for $j=1,2,3$,  $b_{P_2}$ are complex numbers of the form

\begin{equation}\label{bbpp}
 b_{P_2} = \langle G\chi_{\{x:N(x)\in\omega_{P_2}\}}, 
\phi_{P_{1}} \rangle,
\end{equation}
$F_P$ are functions depending on the quartile $P$ and $G$ is an
arbitrary function.

In the treatment of the Walsh bilinear Hilbert transform, we just have
\be{aj-def}
a^{(j)}_{P_j} = \langle f_j, \phi_{P_j} \rangle
\end{equation}

but we will have more sophisticated sequences $a^{(j)}_{P_j}$ 
when dealing with $\Lambda'_{walsh,\P,\Q}$.

In order to estimate these expressions it shall be convenient to
introduce some norms on sequences of tiles (as in
\cite{mtt:walshbiest}, \cite{mtt:fourierbiest}).

\begin{definition}\label{size-def}
Let $\P$ be a finite collection of quartiles, $j=1,2,3$, and let $(a_{P_j})_{P \in \P}$ be a sequence of complex numbers as before.
  We define the \emph{size} of this sequence by
$$ \size_j( (a_{P_j})_{P \in \P} ) := \sup_{T \subset \P}
(\frac{1}{|I_T|} \sum_{P \in T} |a_{P_j}|^2)^{1/2}$$
where $T$ ranges over all trees in $\P$ which are $i$-trees for some $i \neq j$.
We also define the \emph{energy} of the sequence by
$$ \energy_j((a_{P_j})_{P \in \P} ) := \sup_{\D \subset \P}
(\sum_{P \in \D} |a_{P_j}|^2)^{1/2}$$
where $\D$ ranges over all subsets of $\P$ such that the tiles $\{ P_j: P \in \D \}$ are pairwise disjoint.
\end{definition}

The size measures the extent to which the sequence $a_{P_j}$ can concentrate on a single tree and should be thought of as a phase-space variant of the BMO norm.  The energy is a phase-space variant of the $L^2$ norm.  As the  notation suggests, the number $a_{P_j}$ should be thought of as being associated with the tile $P_j$ rather than with the larger quartile $P$.

The usual BMO norm can be written using an $L^2$ oscillation or an $L^1$ oscillation, and the two notions are equivalent thanks to the John-Nirenberg inequality.  The analogous statement for size is (see \cite{mtt:walshbiest}):

\begin{lemma}\label{jn}  Let $\P$ be a finite collection of quartiles, $j=1,2,3$, and let $(a_{P_j})_{P \in \P}$ be a sequence of complex numbers. Then
\be{jn-est}
\size_j( (a_{P_j})_{P \in \P} ) \sim \sup_{T \subset \P}
\frac{1}{|I_T|} \| ( \sum_{P \in T} |a_{P_j}|^2 \frac{\chi_{I_P}}{|I_P|} )^{1/2} \|_{L^{1,\infty}(I_T)}
\end{equation}
where $T$ ranges over all trees in $\P$ which are $i$-trees for some $i \neq j$.
\end{lemma}
We need to define now ``sizes'' and ``energies'' for our 
$b_{P_2}$ sequences. This time, they will no longer depend on the
index ``$j$'' as before.

\begin{definition}\label{bsize-def}
Let $\P$ be a finite collection of quartiles and let 
$(b_{P_2})_{P\in\P}$ be a sequence of complex numbers of the 
form \ref{bbpp}. We define the size of this sequence by

$$\size( (b_{P_2})_{P \in \P} ) :=
\sup_{T\subseteq\P}
\sup_{P\in T}
\sup_{P<P'}
\frac{1}{|I_{P'}|}\int_{I_{P'}}|G|
\chi_{\{x:N(x)\in\omega_{P'_1}\cup\omega_{P'_2}\}}\,dx$$
where $T$ ranges over all trees in $\P$. We also define the energy of this
sequence by

$$\energy( (b_{P_2})_{P \in \P} ) :=
\sup_{\D\subseteq\P}(\sum_{P\in \D}\int_{I_P}
|G|\chi_{\{x:N(x)\in\omega_{P_1}\cup\omega_{P_2}\}})\, dx$$
where $\D$ ranges over all subsets of $\P$ such that the 
corresponding set of sub-bitiles 
$\{I_P\times (\omega_{P_1}\cup\omega_{P_2})\,/\,P\in D\}$
is a set of disjoint bitiles. 
\end{definition}

The following estimate is standard (see \cite{mtt:walshbiest} for a proof)
 and is the main combinatorial tool needed to obtain estimates on 
\eqref{trilinear}.

\begin{proposition}\label{abstract}
Let $\P$ be a finite collection of quartiles, and for each $P \in \P$ and $j=1,2,3$ let $a^{(j)}_{P_j}$ be a complex number as before.  Then
\be{trilinear-est}
| \sum_{P \in \P} \frac{1}{|I_P|^{1/2}} a^{(1)}_{P_1} a^{(2)}_{P_2} a^{(3)}_{P_3}|
\lesssim \prod_{j=1}^3 
\size_j( (a^{(j)}_{P_j})_{P} )^{\theta_j}
\energy_j( (a^{(j)}_{P_j})_{P} )^{1-\theta_j}
\end{equation}
for any $0 \leq \theta_1, \theta_2, \theta_3 < 1$ with $\theta_1 + \theta_2 + \theta_3 = 1$, with the implicit constant depending on the $\theta_j$.
\end{proposition}

If we ignore endpoint issues, Proposition \ref{abstract} says that we can estimate \eqref{trilinear} by taking two of the sequences in the energy norm and the third sequence in the size norm.  This is analogous to the H\"older inequality which asserts that a sum $\sum_i a_i b_i c_i$ can be estimated by taking two sequences in $l^2$ and the third in $l^\infty$.

The main combinatorial tool needed to estimate (\ref{bilinear}) is
the analogue of the above proposition for $b_{P_2}$ sequences:
\begin{proposition}\label{babstract}
Let $\P$ be a finite collection of quartiles, and for each $P\in\P$ 
let $a_{P_1}$ and $b_{P_2}$ complex numbers as before. Then,

\begin{equation}\label{bilinear-est}
|\sum_{P\in\P}a_{P_1}b_{P_2}|\lesssim
\size_1((a_{P_1})_{P})^{\theta_1}
\size((b_{P_2})_{P})^{\theta_2}
\energy_1((a_{P_1})_{P})^{1-\theta_1}
\energy((b_{P_2})_{P})^{1-\theta_2}
\end{equation}
for any $0<\theta_1,\theta_2<1$ with $\theta_1+2 \theta_2=1$.
\end{proposition}
The proof of this Proposition will be presented later on. 
In the meantime, we will take it for granted.
Of course, in order to use Proposition \ref{abstract} we will need some 
estimates on sizes and energies.  In the case when $a^{(j)}_{P_j}$ 
is given by \eqref{aj-def}, the relevant estimates are quite straightforward:
 
\begin{lemma}\label{energy-lemma}
Let $j=1,2,3$, $f_j$ be a function in $L^2(\R)$, and let $\P$ be a finite collection of quartiles.  Then we have
\be{energy-lemma-est}
\energy_j((\langle f_j, \phi_{P_j} \rangle)_{P \in \P} ) \leq
\| f_j \|_2.
\end{equation}
\end{lemma}

\begin{proof}  The wave packets $\phi_{P_j}$ are orthonormal whenever the $P_j$ are disjoint.  The claim then follows immediately from Bessel's inequality.
\end{proof}

\begin{lemma}\label{size-lemma}
Let $j=1,2,3$, $E_j$ be a set of finite measure, $f_j$ be a function in $X(E_j)$, and let $\P$ be a finite collection of quartiles.  Then we have
\be{size-lemma-est}
\size_j( (\langle f_j, \phi_{P_j} \rangle)_{P \in \P} ) \lesssim
\sup_{P \in \P} \frac{|E_j \cap I_P|}{|I_P|}.
\end{equation}
\end{lemma}

\begin{proof}
This shall be a Walsh version of the proof of Lemma 7.8 in \cite{cct}.

From Lemma \ref{jn} it suffices to show that
$$
\| F_T \|_{L^{1,\infty}(I_T)} \lesssim  |I_T|  \sup_{P \in \P} \frac{|E_j \cap I_P|}{|I_P|}
$$
for all $i \neq j$ and all trees $T$, where $F_T$ is the vector-valued function
$$ F_T := (\langle f_j, \phi_{P_j} \rangle \frac{\chi_{I_P}}{|I_P|^{1/2}} )_{P \in T}.$$

It suffices to prove this estimate in the case when $T$ contains its top $P_T$, since in the general case one could then decompose $T$ into disjoint trees with this property and then sum.  In this case it thus suffices to show
$$
\| F_T \|_{L^{1,\infty}(I_T)} \lesssim  |E_j \cap I_T|.$$
From the definition of $F_T$ it is clear that we may restrict $f_j$ and $E_j$ to $I_T$, in which case it suffices to show
$$
\| F_T \|_{L^{1,\infty}(I_T)} \lesssim \|f_j\|_1.$$
We shall assume that $T$ is centered at the frequency origin in the sense that 0 is on the boundary of $w_{P_T}$.  (The general case can then be handled by modulating by an appropriate Walsh ``plane wave'').  But then the linear operator $f_j \mapsto F_T$ is a (vector-valued) dyadic Calder\'on-Zygmund operator, and the claim follows from standard theory.
\end{proof}

Likewise, in order to use Proposition \ref{babstract}, we need again estimates 
on sizes and energies. In the case of the $b_{P_2}$ sequences, these are  also
easy to obtain:

\begin{lemma}\label{benergy-lemma}
Let $f\in L^1(\R)$ and $\P$ be a finite collection of quartiles. 
Then, one has

\begin{equation}\label{benergy-lemma-est}
\energy((
\langle f\chi_{\{x:N(x)\in\omega_{P_2}\}},\phi_{P_1}\rangle)_{P\in\P})
\lesssim \|f\|_1.
\end{equation}
\end{lemma}

\begin{proof}
One just has to observe that the sets 
$\{x\in I_P / N(x)\in\omega_{P_1}\cup\omega_{P_2}\}$
are pairwise disjoint as $P$ varies inside a $D$ as in definition 
\ref{bsize-def}.
\end{proof}

\begin{lemma}\label{bsize-lemma}
Let $E$ be a set of finite measure and $f\in X(E)$. Then,

\begin{equation}\label{bsize-lemma-est}
\size((
\langle f\chi_{\{x:N(x)\in\omega_{P_2}\}},\phi_{P_1}\rangle)_{P\in\P})
\lesssim \sup_{P\in\P}\frac{|E\cap I_P|}{|I_P|}.
\end{equation}
\end{lemma}
The proof follows directly from definitions.

In the next section we shall show how the above size and energy estimates can be combined with  Proposition \ref{babstract} 
and the interpolation theory of the previous section to obtain 
Theorem \ref{carleson-walsh}.  
To prove the estimates for the trilinear operator $T'_{walsh,\P,\Q}$ 
we need some more sophisticated size and energy estimates, which we will pursue after the proof of Theorem \ref{carleson-walsh}.

\section{Proof of Theorem \ref{carleson-walsh}}\label{carleson-walsh-sec}

We now give a proof of Theorem \ref{carleson-walsh}.  We present its proof
here for expository purposes, and also because we shall need 
Theorem \ref{carleson-walsh} to prove the size and energy estimates needed for Theorem \ref{teorema1}.

Fix the collection $\P$ of quartiles, and let $\Lambda$ denote the 
bilinear form
\bas
 \Lambda(f_1,f_2) &:= \langle C_{walsh,\P}(f_1),f_2 \rangle\\
&= \sum_{P \in \P} 
\langle f_1, \phi_{P_1} \rangle
\langle f_2 \chi_{\{x:N(x)\in\omega_{P_2}\}}, \phi_{P_1} \rangle.
\end{align*}
We shall use the notation of Section \ref{interp-sec}, with the obvious modification for bilinear forms as opposed to trilinear forms.

Let us also consider $E_1, E_2$ sets of finite measure and $1<q<2$.
We are going to prove directly that there exists a major subset
$E'_2$ of $E_2$ so that

\begin{equation}\label{1}
|\Lambda(f_1,f_2)|\lesssim |E_1|^{1/q}|E_2|^{1/q'},
\end{equation}
for every $f_1\in X(E_1)$, $f_2\in X(E'_2)$ and also that there exists
a major subset $E'_1$ of $E_1$ so that

\begin{equation}\label{2}
|\Lambda(f_1, f_2)|\lesssim |E_2|,
\end{equation}
for every $f_1\in X(E'_1)$, $f_2\in X(E_2)$\footnote{This is actually
 equivalent to the fact that the adjoint $C^*$ of the linearized Carleson 
operator, is
of weak type $(1,1)$; thus, it differs at this endpoint from 
the Carleson operator, which is not of
weak type $(1,1)$, see \cite{fefferman}.}.

Then, by using the interpolation arguments in \cite{cct}, it follows that
the form $\Lambda$ is of restricted type $\alpha$, for every $\alpha$
in the interior of the segment defined by the endpoints $(0,1)$ and $(1,0)$.
Finally, theorem \ref{carleson-walsh} is implied by the classical 
Marcinkiewicz interpolation
theorem.

It thus remains to prove (\ref{1}) and (\ref{2}).

To prove (\ref{1}), we may assume by scaling invariance, that $|E_2|=1$.
Define the exceptional set

$$\Omega:=\bigcup_{j=1}^2\{x/M\chi_{E_j}>C|E_j|\},$$
where $M$ is the Hardy-Littlewood maximal function. By the classical
Hardy-Littlewood inequality, we have $|\Omega|<1/2$ if $C$ is big enough.
Thus, if we set $E'_2:=E_2\setminus\Omega$, then $E'_2$ is a major subset
of $E_2$. Let now $f_1\in X(E_1)$ and $f_2\in X(E'_2)$. We need to show that

\begin{equation}\label{ab1}
|\sum_{P\in\P}a_{P_1}b_{P_2}|\lesssim |E_1|^{1/q},
\end{equation}
where we denote by

\begin{equation}
\begin{split}
a_{P_1} &:= \langle f_1,\phi_{P_1} \rangle\\
b_{P_2} &:= \langle f_2\chi_{\{x/N(x)\in\omega_{P_2}\}} ,\phi_{P_1} \rangle.
\end{split}
\end{equation}
Also, we may restrict the quartile set $\P$ to those quartiles $P$ for which $I_P \not \subset \Omega$, since $b_{P_2}$ vanishes for all other quartiles.  
By the definition of $\Omega$ we thus have

$$ \frac{|E_1 \cap I_P|}{|I_P|} \lesssim |E_1|$$
for all remaining tiles $P \in \P$.  From Lemma \ref{size-lemma} we thus have

$$ \size_1( (a_{P_1})_{P \in \P} ) \lesssim |E_1|.$$
Also, from Lemma \ref{energy-lemma} and the fact that $f_1 \in X(E'_1)$
 we have
$$ \energy_1( (a_{P_1})_{P \in \P} ) \lesssim |E_1|^{1/2}.$$
Similarly, this time by applying lemma \ref{bsize-lemma} and lemma
\ref{benergy-lemma}, we have

$$\size( (b_{P_2})_{P \in \P} )\lesssim 1$$
and

$$\energy( (b_{P_2})_{P \in \P} )\lesssim 1.$$
From Proposition \ref{babstract} we thus have

$$
|\sum_{P \in \P} 
a_{P_1} b_{P_2}| \lesssim 
|E_1|^{\theta_1}|E_1|^{(1-\theta_1)/2}=|E_1|^{(1+\theta_1)/2},
$$
for every $\theta_1\in (0,1)$. If we chose now $\theta_1$ so that
$(1+\theta_1)/2=1/q$, this proves (\ref{1}).

To prove (\ref{2}), we assume again without loss of generality that
$|E_1|=1$, and define $E'_1$ similarly. This time, we get the bounds

\begin{equation}
\begin{split}
\size_1( (a_{P_1})_{P \in \P} )& \lesssim 1\\
\energy_1( (a_{P_1})_{P \in \P} )& \lesssim 1\\
\size( (b_{P_2})_{P \in \P} )& \lesssim |E_2|\\
\energy( (b_{P_2})_{P \in \P} )& \lesssim |E_2|
\end{split}
\end{equation}
and finally, by applying Proposition \ref{babstract}, we obtain

$$|\sum_{P\in\P}a_{P_1} b_{P_2}|\lesssim |E_2|^{\theta_2}|E_2|^{1-\theta_2}
= |E_2|$$
which completes the proof.

\section{size and energy estimates for $\Lambda'_{walsh,\P,\Q}$}

We now continue the study of the form $\Lambda'_{walsh,\P,\Q}$.  Fix $\P$, $\Q$
and drop any indices $\P$ and $\Q$ for notational convenience.

In the expression $\Lambda'_{walsh}$ the $Q$ tile in the inner summation has a narrower frequency interval, and hence a wider spatial interval, than the $P$ tile in the outer summation.  Thus the inner summation has a poorer spatial localization than the outer sum.  It shall be convenient to reverse the order of summation so that the inner summation is instead more strongly localized spatially than the outer summation.  Specifically, we rewrite $\Lambda'_{walsh}$ as 
$$
\Lambda'_{walsh}(f_1,f_2,f_3,f_4)=
\sum_{Q\in\Q}\frac{1}{|I_Q|^{1/2}} a^{(1)}_{Q_1} a^{(2)}_{Q_2} a^{(3)}_{Q_3}
$$
where
\begin{equation}\label{l'rightform}
\begin{split}
a^{(1)}_{Q_1} &:= \langle f_1,\phi_{Q_1} \rangle\\
a^{(2)}_{Q_2} &:= \langle f_2,\phi_{Q_2} \rangle\\
a^{(3)}_{Q_3} &:= \sum_{P\in\P\,;\,\omega_{Q_3}\subseteq \omega_{P_1}}
\langle f_3\chi_{\{x:N(x)\in\omega_{P_2}\}},\phi_{P_1} \rangle
\langle \phi_{P_1},\phi_{Q_3}\rangle. 
\end{split}
\end{equation}
 
In order to estimate our form $\Lambda'_{walsh}$,
 we need analogues of Lemma \ref{energy-lemma} and Lemma \ref{size-lemma} for $a^{(3)}_{Q_3}$.  The crucial new ingredient in doing so shall be the following simple geometric lemma which allows us to decouple the $P$ and $Q$ variables.

\begin{lemma}\label{biest-trick}Let $i,j=1,2,3$ and let 
$\D$ be a collection of quartiles such that the tiles $\{ Q_j: Q \in \D \}$ are pairwise disjoint.  Let $\P' \subset \P$ denote the set
$$ \P' := \{ P \in \P: P_i \leq Q_j \hbox{ for some } Q \in \D \}.$$
Then for every pair of quartiles $P \in \P$, $Q \in \D$ such that 
$P_i \cap Q_j \neq \emptyset$, we have
$$\omega_{Q_j} \subseteq \omega_{P_i} \hbox{ if and only if } P \in \P'.$$
\end{lemma}

\begin{proof}
Let $P \in \P$, $Q \in \D$ be such that $P_i \cap Q_j \neq \emptyset$.

If $\omega_{Q_j} \subseteq \omega_{P_i}$, then $P_i \leq Q_j$, and so 
$P \in \P'$.  This proves the ``only if'' part.

Now suppose that $\omega_{Q_j} \not \subseteq \omega_{P_i}$, and $Q_j < P_i$.  If $P \in \P'$, then we may find $Q' \in \D$ such that $P_i\leq Q'_j$, hence 
$Q'_j < Q_j$.  But this implies that $Q'_j \cap Q_j \neq \emptyset$, contradicting the disjointness of the $Q_j$.  This proves the ``if'' part.
\end{proof}

The energy estimate is given by the following lemma:

\begin{lemma}\label{carleson-energy}
Let $E_3$ be a set of finite measure and $f_3$ be a function in $X(E_3)$.
  Then we have
\be{energy-carleson-est}
\energy_3((a^{(3)}_{Q_3})_{Q \in \Q}) \lesssim
|E_3|^{1/2}.
\end{equation}
\end{lemma}

\begin{proof}
By Definition \ref{size-def}, we need to show that
\be{ortho-bht}
(\sum_{Q \in \D} |a^{(3)}_{Q_3}|^2)^{1/2} \lesssim 
|E_3|^{1/2}
\end{equation}
for any collection $\D$ of quartiles in $\Q$ such that the tiles 
$\{ Q_3: Q \in \D\}$ are disjoint.  

Fix $\D$, and define the set $\P'$ by
$$ \P' := \{ P \in \P: P_1 \leq Q_3 \hbox{ for some } Q \in \D \}.$$
By Lemma \ref{biest-trick} and \eqref{l'rightform} we may write
$$
a^{(3)}_{Q_3} = \sum_{P\in\P'}
\langle f_3\chi_{\{x:N(x)\in\omega_{P_2}\}},\phi_{P_1} \rangle
\langle \phi_{P_1},\phi_{Q_3}\rangle$$
for all $Q \in \D$.  
We can also rewrite this as
$$ a^{(3)}_{Q_3} = \langle C^*_{walsh,\P'}(f_3), \phi_{Q_3} \rangle$$
where $C^*_{walsh,\P'}$ is  the adjoint of the Carleson operator
 $C_{walsh,\P'}$.  Since the $\phi_{Q_3}$ are orthonormal as 
$Q$ varies in $\D$, we may use Bessel's inequality to estimate the left-hand 
side of \eqref{ortho-bht} by
$$ \| C^*_{walsh,\P'}(f_3) \|_2.$$
The claim then follows from Theorem \ref{carleson-walsh} and the 
assumption $f_3 \in X(E_3)$.
\end{proof}

The analogue of Lemma \ref{size-lemma} is

\begin{lemma}\label{carleson-size}Let $\epsilon>0$,
 $E_3$ be a set of finite measure and $f_3$ be a function in $X(E_3)$.  
Then we have
\be{size-carleson-est}
\size_3((a^{(3)}_{Q_3})_{Q \in \Q}) \lesssim
\sup_{Q \in \Q} (\frac{|E_3 \cap I_Q|}{|I_Q|})^{1/(1+\epsilon)}.
\end{equation}
\end{lemma}

\begin{proof}
By Lemma \ref{jn} it suffices to show that
$$\frac{1}{|I_T|^{1/(1+\epsilon)}}
\| (\sum_{Q \in T} |a^{(3)}_{Q_3}|^2 
\frac{\chi_{I_Q}}{|I_Q|})^{1/2} \|_{L^{1+\epsilon}(I_T)}
\lesssim 
\sup_{Q \in \Q} (\frac{|E_3 \cap I_Q|}{|I_Q|})^{1/(1+\epsilon)}
$$
for some $i \neq 3$ and some $i$-tree $T$.  
One may also reduce the above inequality to:
\be{equation}\label{weak'}
\| (\sum_{Q \in T} |a^{(3)}_{Q_3}|^2 \frac{\chi_{I_Q}}{|I_Q|})^{1/2} 
\|_{L^{1+\epsilon}(I_T)}
\lesssim 
|E_3 \cap I_T|^{1/(1+\epsilon)}.
\end{equation}
From \eqref{l'rightform} we see that the only quartiles $P \in \P$ which matter are those such that $I_P \subseteq I_T$.  
Thus we may restrict $f_3$, $E_3$, to $I_T$.

Fix $i$, $T$, and define the set $\P'$ by
$$ \P' := \{ P \in \P: P_1 \leq Q_3 \hbox{ for some } Q \in T \}.$$
By Lemma \ref{biest-trick} and \eqref{l'rightform} as before we have
$$ a^{(3)}_{Q_3} = \langle C^*_{walsh,\P'}(f_3), \phi_{Q_3} \rangle$$
for all $Q \in T$.  By the dyadic Littlewood-Paley estimate for the tree $T$ 
we may thus reduce \eqref{weak'} to
\[
\| C^*_{walsh,\P'}(f_3) \|_{L^{1+\epsilon}(I_T)}
\lesssim 
|E_3|^{1/(1+\epsilon)}.
\]
But this follows from Theorem \ref{carleson-walsh} and the 
assumption $f_3 \in X(E_3)$.
\end{proof}

\section{proof of theorem \ref{teorema1}}

We can now present the proof of Theorem \ref{teorema1}.
Fix the collections $\P$ and $\Q$ of quartiles.
We first show that $\Lambda'_{walsh,\P,\Q}$
is of restricted weak type $\alpha$ for all admissible 3-tuples 
$(\alpha_1, \alpha_2, \alpha_3)$ 
 arbitrarily close to $A_2$, $A_3$, so that the bad index is 3.

Fix $\alpha$ as above and let $E_1$, $E_2$, $E_3$ be sets of finite measure.  

By scaling invariance we may assume that $|E_3|=1$. 
We need to find a major subset $E'_3$ of $E_3$ such that
$$
|\Lambda'_{walsh,\P,\Q}(f_1, f_2, f_3)|\lesssim |E|^{\alpha}
$$
for all functions $f_i\in X(E'_i)$, $i=1,2,3$.

Define the exceptional set $\Omega$ by
\[\Omega := \bigcup_{j=1}^{3}\{M\chi_{E_j}>C |E_j|\}\]
where $M$ is the dyadic Hardy-Littlewood maximal  function.
By the classical Hardy-Littlewood inequality, we have $|\Omega|<1/2 $
if $C$ is a sufficiently large constant.  Thus if we set 
$E'_3 := E_3 \setminus \Omega$, then $E'_3$ is a major subset of $E_3$.

Let than $f_i \in X(E'_i)$ for $i=1,2,3$.  We need to show
$$
|\sum_{Q \in \Q} \frac{1}{|I_Q|^{1/2}}
a^{(1)}_{Q_1} a^{(2)}_{Q_2} a^{(3)}_{Q_3}| \lesssim |E|^\alpha$$
where $a^{(j)}_{Q_j}$ is defined by \eqref{l'rightform}.

We may restrict the quartile set $\Q$ to those quartiles $Q$ for which 
$I_Q \not \subset \Omega$, since $a^{(3)}_{Q_j}$ vanishes for all other quartiles.  By the definition of $\Omega$ we thus have

$$ \frac{|E_j \cap I_Q|}{|I_P|} \lesssim |E_j|$$
for all remaining tiles $Q \in \Q$ and $j=1,2,3$.  
From Lemma \ref{size-lemma} we thus have
$$ \size_j( (a^{(j)}_{Q_j})_{Q \in \Q} ) \lesssim |E_j|$$
for $j=1,2$.  
On the other hand, by Lemma \ref{carleson-size} one has
$$ \size_3( (a^{(3)}_{Q_3})_{Q \in \Q} ) \lesssim 1$$
since $|E_3|=1$.
Also, from Lemma \ref{energy-lemma}, Lemma \ref{carleson-energy}
and the fact that $f_j \in X(E'_j)$ we have
$$ \energy_j( (a^{(j)}_{Q_j})_{Q \in \Q} ) \lesssim |E_j|^{1/2},$$
for $j=1,2,3.$
From Proposition \ref{abstract} we thus have
$$
|\sum_{Q \in \Q} \frac{1}{|I_Q|^{1/2}}
a^{(1)}_{Q_1} a^{(2)}_{Q_2} a^{(3)}_{Q_3}| \lesssim 
\prod_{j=1}^2 |E_j|^{(1-\theta_j)/2} |E_j|^{\theta_j} $$
for any $0 \leq \theta_1, \theta_2< 1$ such that there exists
$0\leq \theta_3<1$ with $\theta_1 + \theta_2 + \theta_3 = 1$.  
The claim then follows by choosing $\theta_1 := 2\alpha_1-1$, $\theta_2 := 2\alpha_2-1$ ; note that there exist choices of $\alpha$ arbitrarily close to 
$A_2$ or $A_3$, for which the constraints 
on $\theta_1, \theta_2, \theta_3$ are satisfied.  

To prove the restricted type estimates for $\alpha$ arbitrarily
close to $A_4, A_5, A_6, A_1$, one argues in the same way, by taking
advantage of the fact that $\epsilon$ in Lemma \ref{carleson-size}
can be arbitrarily small.

This concludes the proof of Theorem \ref{teorema1}.

\section{Tile norms, part 2}

In this section we begin the study of $\Lambda''_{walsh,\P,\Q}$.
As before, fix $\P,\Q$ and drop any indices $\P$ and $\Q$ for
notational convenience. Also, as in the previous sections, it is more 
convenient to rewrite $\Lambda''_{walsh}$ as

\[\Lambda''_{walsh}(f_1,f_2,f_3)=\]
\begin{equation}\label{Lambda''walsh}
\sum_{Q\in\Q}
\langle f_2, \phi_{Q_1}\rangle
\left< \phi_{Q_1} \chi_{\{x:N(x)\in \omega_{Q_2}\}}
,\sum_{P\in\P; \omega_{Q_1}\subseteq\omega_{P_2}}
\langle f_1, \phi_{P_1}\rangle\phi_{P_1}
\chi_{\{x:N(x)\in \omega_{P_2}\}}f_3\right>.
\end{equation}
Unfortunately, this form $\Lambda''_{walsh}$ is not as symmetric as
$\Lambda'_{walsh}$ and to deal with it will require some notational
adjustments.
First, for $j=1,2,3$ and $F$ arbitrary, we will write for simplicity
$$\size_{j,\P}(F):=\size_j(\langle F, \phi_{P_j}\rangle_P)$$
or sometimes
$$\size_{j,\Q}(F):=\size_j(\langle F, \phi_{Q_j}\rangle_Q),$$
depending which family of functions between
$(\phi_{P_j})_P$ and $(\phi_{Q_j})_Q$ is the relevant one. In the same
way we will write $\energy_{j,\P}(F)$ for 
$\energy_j(\langle F, \phi_{P_j}\rangle_P)$ and $\energy_{j,\Q}(F)$ for
$\energy_j(\langle F, \phi_{Q_j}\rangle_Q)$.

Similarly, for an arbitrary $G$ we will use the notation
$$\size_{\P}(G):=\size(\langle G\chi_{\{x:N(x)\in \omega_{P_2}\}},\phi_{P_1} 
\rangle_P)$$
or sometimes
$$\size_{\Q}(G):=\size(\langle G\chi_{\{x:N(x)\in \omega_{Q_2}\}},\phi_{Q_1} 
\rangle_Q)$$
and also $\energy_{\P}(G)$ for
$\energy(\langle G\chi_{\{x:N(x)\in \omega_{P_2}\}},\phi_{P_1} 
\rangle_P)$.
 For simplicity, we will usually omit the dependence on $\P$ and $\Q$ 
when this will
be clear from the context.

 We need two more definitions of ``sizes'' and
``energies''
in order to handle the form $\Lambda''_{walsh}$.
In fact, these definitions have been inspired by the
right hand side term in (\ref{Lambda''walsh}).

\begin{definition}\label{size'energy'-def}
Let $\P$ and $\Q$ be two finite collections of quartiles and let
$f_1$ and $f_3$ be two arbitrary functions. We define the $\size'$ of the
pair $(f_1,f_3)$ by
\begin{equation}\label{size'-def}
\size'(f_1,f_3)=\size'_{\Q,\P}(f_1,f_3):=
\end{equation}
$$
\sup_{T\subseteq\Q}
\sup_{Q\in T}
\sup_{Q':Q_{12}<Q'_{12}}
\frac{1}{|I_{Q'}|}
\int_{I_{Q'}}
|C^c_{walsh,Q',\P}(f_1)f_3|
\chi_{\{x:N(x)\in \omega_{Q'_1}\cup\omega_{Q'_2}\}}\,dx\footnote{This is like 
a phase plane version of $\|C(f_1)f_3\|_{\infty}$.}
$$
where
\begin{equation}\label{carlesonc}
C^c_{walsh,Q',\P}(f_1):=
\sum_{P\in\P;\omega_{Q'_1}\subseteq\omega_{P_2}}
\langle f_1, \phi_{P_1}\rangle\phi_{P_1}
\chi_{\{x:N(x)\in \omega_{P_2}\}}
\end{equation}
and $T$ ranges over all trees in $\Q$.
We also define the $\energy'$
of the pair $(f_1,f_3)$ by
\begin{equation}\label{energy'-def}
\energy'(f_1,f_3)=\energy_{\Q,\P}'(f_1,f_3):=
\end{equation}
$$
\sup_{\D\subseteq \Q}
\left(\sum_{Q\in\D}
\int_{I_{Q}}
|C^c_{walsh,Q,\P}(f_1)f_3|
\chi_{\{x:N(x)\in \omega_{Q_1}\cup\omega_{Q_2}\}}\,dx\right)\footnote{This is 
like a phase plane version of $\|C(f_1)f_3\|_1$.}
$$
where $\D$ ranges over all subsets of $\Q$ such that the corresponding
set of sub-bitiles $Q_{12}$
is a set of disjoint bitiles.
\end{definition}

\begin{definition}\label{size''energy''-def}
Let $\P$ and $\Q$ be two finite collections of quartiles and let
$f_3$ and $f_1$ be two arbitrary functions. We define the $\size''_1$ of the
pair $(f_3,f_1)$ by
\begin{equation}\label{size''-def}
\size''_1(f_3,f_1)=\size''_{1,\Q,\P}(f_3,f_1):=
\end{equation}
$$
\sup_{T\subseteq\Q}
\frac{1}{|I_T|}
\|
\sum_{Q\in T}
\langle f_3,\phi_{Q_1}\rangle\phi_{Q_1}
\sum_{P\in\tilde{T};\omega_{Q_1}\subseteq\omega_{P_2}}
\langle f_1,\phi_{P_1}\rangle\phi_{P_1}
\|_1
$$
where $T$ ranges over all $i$-trees in $\Q$, $i\neq 1$ and 
$\tilde{T}$ is the set of all $P\in\P$ for which there exist $Q'$ and
$Q''$ in $T$ so that $Q'_2\leq P_2\leq Q''_2$.
We also define the $\energy''_1$
of the pair $(f_3,f_1)$ by

\begin{equation}\label{energy''-def}
\energy''_1(f_3,f_1)=\energy''_{1,\Q,\P}(f_3,f_1):=
\end{equation}
$$
\sup_{\D\subseteq \Q}
\left(\sum_{T\in\D}
\|
\sum_{Q\in T}
\langle f_3,\phi_{Q_1}\rangle\phi_{Q_1}
\sum_{P\in\tilde{T};\omega_{Q_1}\subseteq\omega_{P_2}}
\langle f_1,\phi_{P_1}\rangle\phi_{P_1}
\|^2_2\right)^{1/2}
$$
where $\D$ ranges over all subsets of $\Q$ 
which are unions of $i$-trees $i\neq 1$ so that for every
$T,T'\in D$, $T$ and $T'$ are disjoint and also $\tilde{T}$ and
$\tilde{T'}$ are disjoint.
\end{definition}

The following estimate is the main combinatorial tool needed
to estimate our form $\Lambda''_{walsh}$ (\ref{Lambda''walsh}).

\begin{proposition}\label{moreabstract}
Let $\P$ and $\Q$ be finite collections of quartiles and $f_1,f_2,f_3$
be three arbitrary functions so that $\|f_j\|_{\infty}\lesssim 1$
for $j=1,2,3$. Then,
\bas
|\Lambda''_{walsh}(f_1,f_2,f_3)|& \lesssim
\size_1(f_2)^{\alpha_1}\size'(f_1,f_3)^{\alpha_2}
\energy_1(f_2)^{1-\alpha_1}\energy'(f_1,f_3)^{1-\alpha_2}\\
&+\min(E_1(f_1,f_2,f_3,(\beta_j)_j), E_2(f_1,f_2,f_3,(\gamma_j)_j))
\end{align*}
where
$$E_1(f_1,f_2,f_3,(\beta_j)_j)=\size''_1(f_2,f_1)^{\beta_1}
\size(f_3)^{\beta_2}\energy''_1(f_2,f_1)^{1-\beta_1}
\energy(f_3)^{1-\beta_2}$$
and
$$E_2(f_1,f_2,f_3,(\gamma_j)_j)=\size_1(f_2)^{\gamma_1}
\size(f_3)^{\gamma_3}\energy_1(f_2)^{1-\gamma_1}
\energy_1(f_1)^{1-\gamma_2}\energy(f_3)^{1-\gamma_3},$$
for any $\alpha_j,\beta_j,\gamma_j\in (0,1)$ so that
$\alpha_1+2\alpha_2=1$, $\beta_1+2\beta_2=1$ and 
$\gamma_1+\gamma_2+2\gamma_3=2$.
\end{proposition}

As before, for the moment we take this Proposition for granted.
To use it,
one needs more size and energy estimates which we will present in the next
section.

\section{Size and energy estimates for $\Lambda''_{walsh,\P,\Q}$}

The first lemma will be useful for proving Theorem \ref{teorema2} 
near the vertices $M_{56}$, $M_{12}$:

\begin{lemma}\label{carleson'-energy}
Let $E_j$ be sets of finite measure and $f_j$ be functions in $X(E_j)$ for 
$j=1,3$.  Then we have
\be{energy-carleson'-est}
\energy'(f_1,f_3) \lesssim
|E_1|^{1-\theta} |E_3|^{\theta}
\end{equation}
for any $0 < \theta < 1$, with the implicit constant depending on $\theta$.
\end{lemma}

\begin{proof}
By Definition \ref{size'energy'-def}, we need to show that
\be{ortho-carleson'}
\sum_{Q\in \D}
\int_{I_Q}|C^c_{walsh,Q,\P}(f_1)f_3|
\chi_{\{x:N(x)\in\omega_{Q_1}\cup\omega_{Q_2}\}}\,dx
\lesssim |E_1|^{1-\theta} |E_3|^{\theta}
\end{equation}
for any collection $\D$ of quartiles as in Definition 
\ref{size'energy'-def}.

Fix $\D$, and define the set $\P'$ by
$$ \P' := \{ P \in \P: P_2 \leq Q_1\hbox{ for some } Q \in \D \}.$$
Since the sets 
$\{x\in I_Q/N(x)\in\omega_{Q_1}\cup\omega_{Q_2}\}$ are pairwise disjoint as
$Q$ varies inside $\D$, and by using Lemma \ref{biest-trick}, one can
majorize the left hand side of (\ref{ortho-carleson'}) by

$$\|C_{walsh,\P'}(f_1)f_3\|_1.$$
The claim then follows from H\"{o}lder inequalities, 
Theorem \ref{carleson-walsh} and the assumptions
$f_1 \in X(E_1)$, $f_3 \in X(E_3)$.
\end{proof}

To prove Theorem \ref{teorema2} near $A_2$ 
we will use the following variant:

\begin{lemma}\label{carleson''-energy}
Let $E_j$ be sets of finite measure and $f_j$ be functions in $X(E_j)$ for 
$j=1,3$.  Then we have
\be{energy-carleson''-est}
\energy'(f_1,f_3) \lesssim
(|E_1|\sup_{P \in \P} \frac{|E_3 \cap I_P|}{|I_P|})^{\theta}
|E_3|^{1-\theta}
(\sup_{P \in \P} \frac{|E_1 \cap I_P|}{|I_P|})^{1-2\theta}
\end{equation}
for any $0 < \theta < 1/2$, with the implicit constant depending on $\theta$.
\end{lemma}

\begin{proof}
By repeating the proof of Lemma \ref{carleson'-energy}, we reduce to 
showing that

$$\|C_{walsh,\P'}(f_1)f_3\|_1\lesssim
(|E_1|\sup_{P \in \P} \frac{|E_3 \cap I_P|}{|I_P|})^{\theta}
|E_3|^{1-\theta}
(\sup_{P \in \P} \frac{|E_1 \cap I_P|}{|I_P|})^{1-2\theta}$$
where $\P'$ is an arbitrary subset of $\P$.  
By duality we may write the left-hand side as
$$
|\sum_{P \in \P'} 
\langle f_1, \phi_{P_1} \rangle
\langle f_3 F\chi_{\{x:N(x)\in\omega_{P_2}\}}, \phi_{P_1} \rangle|
$$
for some $L^{\infty}$-normalized function $F$.  
By Proposition \ref{babstract} we may estimate this by
\bas
&\size_1( (\langle f_1, \phi_{P_1} \rangle)_{P \in \P'}))^{\theta_1}\\
&\size_2( (\langle f_3 F\chi_{\{x:N(x)\in\omega_{P_2}\}} , \phi_{P_1}
\rangle)_{P \in \P'}))^{\theta_2}\\
& \energy_1( (\langle f_1, \phi_{P_1} \rangle)_{P \in \P'})^{1-\theta_1}\\
&(\energy_2( (\langle f_3 F\chi_{\{x:N(x)\in\omega_{P_2}\}}, 
\phi_{P_1}\rangle)_{P \in \P'})^{1-\theta_2}
\end{align*}
whenever $0<\theta_1, \theta_2 <1$ and $\theta_1+2\theta_2=1$.
The claim then follows from Lemma \ref{size-lemma}, Lemma \ref{energy-lemma},
 Lemma \ref{bsize-lemma} and  Lemma \ref{benergy-lemma}.
\end{proof}

We also need $\size'$ bounds.

\begin{lemma}\label{carleson'-size}
Let $E_j$ be sets of finite measure and $f_j$ be functions in $X(E_j)$ for 
$j=1,3$.  Then we have
\be{size-carleson'-est}
\size'(f_1,f_3) \lesssim
\sup_{Q \in \Q} (\frac{|E_1 \cap I_Q|}{|I_Q|})^{1-\theta} 
(\frac{|E_3 \cap I_Q|}{|I_Q|})^\theta
\end{equation}
for any $0 < \theta < 1$, with the implicit constant depending on $\theta$.
\end{lemma}

\begin{proof}
By Definition \ref{size'energy'-def} we need to show that for any
$Q'\in\Q$ one has
\begin{equation}\label{pac}
\frac{1}{|I_{Q'}|}
\int_{I_{Q'}}|C^c_{walsh,Q',\P}(f_1)f_3|
\chi_{\{x:N(x)\in\omega_{Q'_1}\cup\omega_{Q'_2}\}}\,dx
\lesssim
\sup_{Q \in \Q} (\frac{|E_1 \cap I_Q|}{|I_Q|})^{1-\theta} 
(\frac{|E_3 \cap I_Q|}{|I_Q|})^\theta.
\end{equation}
Fix $Q'\in\Q$ and define as before the set $\P'$ by
$$\P'=\{P\in\P: P_2\leq Q'_1\}.$$
Arguing as in the proofs of Lemmas \ref{carleson'-energy},
\ref{carleson-size} one can reduce (\ref{pac}) to proving
$$\|C_{walsh,\P'}(f_1)f_3\|_1\lesssim
|E_1|^{1-\theta}|E_3|^{\theta}.$$
Then again the claim follows from H\"{o}lder, Theorem
\ref{carleson-walsh}
and the assumption $f_j\in X(E_j)$.
\end{proof}
Finally, we will prove $\size''_1$ and $\energy''_1$ bounds.

\begin{lemma}\label{carleson'''-size}
Let $E_j$ be sets of finite measure and $f_j$ be functions in $X(E_j)$ for 
$j=1,2$.  Then we have
\begin{equation}
\size''_1(f_2,f_1) \lesssim
\sup_{Q \in \Q} (\frac{|E_2 \cap I_Q|}{|I_Q|})^{1-\theta} 
(\frac{|E_1 \cap I_Q|}{|I_Q|})^\theta
\end{equation}
for any $0 < \theta < 1$, with the implicit constant depending on $\theta$.
\end{lemma}

\begin{proof}
Let $T$ be an $i$ tree in $\Q$, $i\neq 1$. We need to show that
\begin{equation}\label{paraprodus}
\frac{1}{|I_T|}
\|
\sum_{Q\in T}
\langle f_2,\phi_{Q_1}\rangle\phi_{Q_1}
\sum_{P\in\tilde{T};\omega_{Q_1}\subseteq\omega_{P_2}}
\langle f_1,\phi_{P_1}\rangle\phi_{P_1}
\|_1
\lesssim
\sup_{Q \in \Q} (\frac{|E_2 \cap I_Q|}{|I_Q|})^{1-\theta} 
(\frac{|E_1 \cap I_Q|}{|I_Q|})^\theta.
\end{equation}
The left hand side of (\ref{paraprodus}) can be written as
$$
\frac{1}{|I_T|}
\|\Pi_T(
\sum_{Q\in T}
\langle f_2,\phi_{Q_1}\rangle\phi_{Q_1},
\sum_{P\in\tilde{T}}
\langle f_1,\phi_{P_1}\rangle\phi_{P_1})
\|_1$$
where $\Pi_T$ is the Walsh paraproduct associated to $T$ 
(for more details about the classical theory of paraproducts see \cite{stein},
or \cite{auscher} for their Walsh analogues used in this paper).
Since $\Pi_T$ maps $L^p\times L^q\rightarrow L^1$ as long as
$1<p,q<\infty$, $1/p+1/q=1$ (see \cite{stein}, \cite{auscher}), 
the claim follows
by using the fact that $f_j\in X(E_j)$ for $j=1,2$.
\end{proof}

\begin{lemma}\label{carleson'''-energy}
Let $E_j$ be sets of finite measure and $f_j$ be functions in $X(E_j)$ for 
$j=1,2$.  Then we have
\begin{equation}
\energy''_1(f_2,f_1)\lesssim
\sup_{Q \in \Q} |E_2|^{(1-\theta)/2}|E_1|^{\theta/2}
\end{equation}
for any $0 < \theta < 1$, with the implicit constant depending on $\theta$.
\end{lemma}

\begin{proof}
Let $\D\subseteq\Q$ be a set as in Definition
\ref{size''energy''-def}.
Using the notation in the proof of the previous lemma, we can write

\begin{equation}\label{paraprodus'}
\energy''_1(f_2,f_1)^2\lesssim \sum_{T\in\D}
\|\Pi_T(
\sum_{Q\in T}
\langle f_2,\phi_{Q_1}\rangle\phi_{Q_1},
\sum_{P\in\tilde{T}}
\langle f_1,\phi_{P_1}\rangle\phi_{P_1})
\|_2^2.
\end{equation}
On the other hand, the right hand side term in 
(\ref{paraprodus'}) can be estimated either by
$$\sum_{T\in\D}
\|
\sum_{Q\in T}
\langle f_2,\phi_{Q_1}\rangle\phi_{Q_1}\|_2^2
\|\sum_{P\in\tilde{T}}
\langle f_1,\phi_{P_1}\rangle\phi_{P_1}
\|_{BMO}^2$$
or by

$$
\sum_{T\in\D}
\|
\sum_{Q\in T}
\langle f_2,\phi_{Q_1}\rangle\phi_{Q_1}\|_{BMO}^2
\|\sum_{P\in\tilde{T}}
\langle f_1,\phi_{P_1}\rangle\phi_{P_1})
\|_2^2$$
by using the well known estimates on paraproducts (see \cite{stein}).
We also observe that the maps 
$f_2\rightarrow \sum_{Q\in T}\langle f_2,\phi_{Q_1}\rangle\phi_{Q_1}$
and
$f_1\rightarrow \sum_{P\in\tilde{T}}\langle
f_1,\phi_{P_1}\rangle\phi_{P_1}$
are both discrete versions of the Hilbert transform and therefore,
they are bounded from $L^{\infty}$ into $BMO$.

Consequently, we obtain two bounds for $\energy''_1(f_2,f_1)$, namely
$$\energy''_1(f_2,f_1)^2\lesssim |E_2|,\,|E_1|,$$
by taking into account the fact that $f_j\in X(E_j)$ and by Bessel's
inequality. The proof ends by interpolating between the above two
estimates.
\end{proof}

\section{Proof of Theorem \ref{teorema2} for the vertex $M_{56}$}

Let $\alpha=(\alpha_1,\alpha_2,\alpha_3)$ be a good admissible tuple
in the interior of the segment $(M_{56}, M_{34})$, very close
to $M_{56}$ (in particular, $\alpha_2=0$, $\alpha_1+\alpha_3=1$ and
$\alpha_1$ is small).

Let us also
fix $E_1,E_2,E_3$ arbitrary sets of finite measure and assume without loss
of generality that $|E_2|=1$.

As before, we define 
\[\Omega := \bigcup_{j=1}^{3}\{M\chi_{E_j}>C|E_j|\}\]
for a large constant $C$, and set $E'_2 := E_2 \setminus \Omega$.  We now fix $f_i \in X(E'_i)$ for $i=1,2,3$.  Our task is then to show

\begin{equation}\label{inegalitatea'}
|\Lambda''_{walsh}(f_1,f_2,f_3)| 
|\lesssim |E|^\alpha.
\end{equation}

As before, we may restrict the collection $\Q$ to those quartiles $Q$
for which $I_Q \not \subset \Omega$, since our sum
vanishes for all other quartiles.
\footnote{Note however that we cannot restrict $\P$ this way, as 
$I_P\subset I_Q$ and $I_P \subset \Omega$
does not imply $I_Q \subset \Omega$.}  This implies that
$$ \frac{|E_j \cap I_Q|}{|I_Q|} \lesssim |E_j|$$
for all remaining tiles $Q \in \Q$ and $j=1,2,3$.  
From these inequalities and from Lemma \ref{size-lemma}, 
Lemma \ref{energy-lemma}, 
Lemma \ref{carleson'-size} and Lemma \ref{carleson'-energy}
 we thus have
\bas
\size_1(f_2) &\lesssim 1\\
\energy_1(f_2) &\lesssim 1\\
\size'(f_1,f_3) &\lesssim |E_1|^{1-\theta}|E_3|^{\theta}\\
\energy'(f_1,f_3) &\lesssim |E_1|^{1-\theta}|E_3|^{\theta}
\end{align*}
for some $0 < \theta < 1$ which we will choose later.  

On the other hand, from the same inequalities and by using
Lemma \ref{bsize-lemma}, Lemma \ref{benergy-lemma}, Lemma
\ref{carleson'''-size} and Lemma \ref{carleson'''-energy}
we have
\bas
\size (f_3) &\lesssim 1\\
\energy (f_3) &\lesssim |E_3|\\
\size''(f_2,f_1) &\lesssim |E_2|^{1-\theta'}|E_1|^{\theta'}=|E_1|^{\theta'}\\
\energy''(f_2,f_1) &\lesssim |E_2|^{(1-\theta')/2}|E_1|^{\theta'/2}=
|E_1|^{\theta'/2}
\end{align*}
where again $\theta'\in (0,1)$ will be chosen later.

By Proposition \ref{moreabstract} we thus can bound the left-hand side of 
\eqref{inegalitatea'} by

$$|E_1|^{1-\theta}|E_3|^{\theta}+|E_1|^{(\beta_1/2+1/2)\theta'}
|E_3|^{1-\beta_2}$$
where $\beta_1,\beta_2\in (0,1)$ are so that $\beta_1+2\beta_2=1$.
Then, one defines $1-\beta_2=\theta=\alpha_3$ and choses $\theta'$
such that $(\beta_1/2+1/2)\theta'+1-\beta_2=1$. Note that all the
constraints on $\theta,\theta',\beta_1,\beta_2$ are satisfied.
After that, we obtain the majorant
$$|E_1|^{\alpha_1}|E_3|^{\alpha_3}$$
and this finishes the proof.

\section{Proof of Theorem \ref{teorema2} for the vertex $M_{12}$}

Let $\alpha=(\alpha_1,\alpha_2,\alpha_3)$ be a good admissible tuple,
very close to $M_{12}$.

Let us also
fix $E_1,E_2,E_3$ arbitrary sets of finite measure and assume without loss
of generality that $|E_3|=1$.

As usual, we define 
\[\Omega := \bigcup_{j=1}^{3}\{M\chi_{E_j}>C|E_j|\}\]
for a large constant $C$, and set $E'_3 := E_3 \setminus \Omega$.  
We now fix $f_i \in X(E'_i)$ for $i=1,2,3$.  Our task is then to show

\begin{equation}\label{inegalitatea''}
|\Lambda''_{walsh}(f_1,f_2,f_3)| 
|\lesssim |E|^\alpha.
\end{equation}

We may also restrict the collection $\P$ to those quartiles $P$
for which $I_P \not \subset \Omega$, since our sum
vanishes for all other quartiles. Since a pair of tiles $(P,Q)$ gives
a nonzero term in our sum iff $I_P\subseteq I_Q$, we may also restrict
$\Q$ to those quartiles $Q$ for which $I_Q\not\subset\Omega$.
  As before, this implies that
$$ \frac{|E_j \cap I_Q|}{|I_Q|} \lesssim |E_j|$$
for all remaining tiles $Q \in \Q$ and $j=1,2,3$.  
From these inequalities and from Lemma \ref{size-lemma}, 
Lemma \ref{energy-lemma}, 
Lemma \ref{carleson'-size} and Lemma \ref{carleson'-energy}
 we thus have
\bas
\size_1(f_2) &\lesssim |E_2|\\
\energy_1(f_2) &\lesssim |E_2|^{1/2}\\
\size'(f_1,f_3) &\lesssim |E_1|^{1-\theta}\\
\energy'(f_1,f_3) &\lesssim |E_1|^{1-\theta}
\end{align*}
for some $0 < \theta < 1$ which we will choose later.  

On the other hand, from the same inequalities and by using
Lemma \ref{bsize-lemma}, Lemma \ref{benergy-lemma}, Lemma
\ref{carleson'''-size} and Lemma \ref{carleson'''-energy}
we have
\bas
\size (f_3) &\lesssim 1\\
\energy (f_3) &\lesssim 1\\
\size''(f_2,f_1) &\lesssim |E_2|^{1-\theta'}|E_1|^{\theta'}=|E_1|^{\theta'}\\
\energy''(f_2,f_1) &\lesssim |E_2|^{(1-\theta')/2}|E_1|^{\theta'/2}=
|E_1|^{\theta'/2}
\end{align*}
where again $\theta'\in (0,1)$ will be chosen later.

By applying Proposition \ref{moreabstract} 
we thus can bound the left-hand side of 
\eqref{inegalitatea''} by

$$|E_2|^{\delta/2+1/2}|E_1|^{1-\theta}+
|E_2|^{(1-\theta')\beta_1+(1-\theta')(1-\beta_1)/2}
|E_1|^{\theta'\beta1+\theta'(1-\beta_1)/2}$$
where $\beta_1,\beta_2\in (0,1)$ are so that $\beta_1+2\beta_2=1$ and 
$\delta\in (0,1)$ and the claim follows by
setting $\delta/2+1/2=
(1-\theta')\beta_1+(1-\theta')(1-\beta_1)/2
=\alpha_2$ and 
$1-\theta=\theta'\beta_1+\theta'(1-\beta_1)/2=\alpha_1$.
The reader may verify that the constraints on $\theta,\theta',\delta,
\beta_1,\beta_2$ can be obeyed for $\alpha$ arbitrarily close
to $M_{12}$.

\section{Proof of Theorem \ref{teorema2} for the vertex $A_2$}

Let $\alpha=(\alpha_1,\alpha_2,\alpha_3)$ be an admissible tuple
very close to the point $A_2$ (in particular, the bad index is $3$).

Let us also
fix $E_1,E_2,E_3$ arbitrary sets of finite measure and assume once again
without loss
of generality that $|E_3|=1$.

As before, we define 
\[\Omega := \bigcup_{j=1}^{3}\{M\chi_{E_j}>C|E_j|\}\]
for a large constant $C$, and set $E'_3 := E_3 \setminus \Omega$.  
We now fix $f_i \in X(E'_i)$ for $i=1,2,3$.  Our task is then to show

\begin{equation}\label{inegalitatea'''}
|\Lambda''_{walsh}(f_1,f_2,f_3)
|\lesssim |E|^\alpha.
\end{equation}

As before, we may restrict our collections
$\P$ and $\Q$ to those quartiles $P$ and $Q$ for which 
$I_P \not \subset \Omega$ and
$I_Q \not \subset \Omega$, since the corresponding terms in our sum 
vanish for all other quartiles.  
This implies in particular that
$$ \frac{|E_j \cap I_P|}{|I_P|} \lesssim |E_j|$$
for all remaining tiles $P \in \P$ and $j=1,2,3$ 
and similarly
$$ \frac{|E_j \cap I_Q|}{|I_Q|} \lesssim |E_j|$$
for all remaining tiles $Q \in \Q$ and $j=1,2,3$.  
From these inequalities and form 
Lemma \ref{size-lemma}, Lemma \ref{energy-lemma}, 
Lemma \ref{carleson'-size} and Lemma \ref{carleson''-energy}
 we  have
\bas
\size_1( f_2) &\lesssim |E_2|\\
\energy_1( f_2 ) &\lesssim |E_2|^{1/2}\\
\size_2(f_1,f_3) &\lesssim |E_1|^{1-\theta}\\
\energy_2(f_1,f_3) &\lesssim |E_1|^{1-\theta}
\end{align*}
for some $0 < \theta < 1/2$ which we will choose later 
(we also used the fact that $|E_3|=1$).  

On the other hand, by using again the above inequalities and also
Lemma \ref{size-lemma}, Lemma \ref{energy-lemma}, Lemma
\ref{bsize-lemma} and Lemma \ref{benergy-lemma} we have
\bas
\size_1(f_2)&\lesssim |E_2|\\
\energy_1(f_2)&\lesssim |E_2|^{1/2}\\
\energy_1(f_1)&\lesssim |E_1|^{1/2}\\
\size(f_3)&\lesssim 1\\
\energy(f_3)&\lesssim 1.
\end{align*}

By Proposition \ref{moreabstract} we can thus bound the left-hand side of 
\eqref{inegalitatea'''} by

\begin{equation}
|E_1|^{1-\theta}|E_2|^{\delta/2+1/2}+
|E_1|^{1/2-\gamma_2/2}|E_2|^{1/2+\gamma_1/2},
\end{equation}
where $0<\gamma_1,\gamma_2<1$. Now our claim follows by setting
$1-\theta=1/2-\gamma_2/2=\alpha_1$ and $\delta/2+1/2=1/2+\gamma_1/2
=\alpha_2$. The reader may check again that the constraints on
$\theta,\gamma_1,\gamma_2$ are satisfied for $\alpha$ arbitrarily
close to $A_2$. This ends the proof.

\section{Combinatorial Lemmas}\label{abstract-sec}

In order to prove Propositions \ref{babstract} and \ref{moreabstract},
we first need to prove certain combinatorial lemmas.

Fix the collections $\P$ and $\Q$.

  We begin by considering the contribution of a single tree:

\begin{lemma}[Tree estimate]\label{single-tree}
Let $T$ be a tree in $\Q$, and $f_1,f_2,f_3$ be three functions as before.  

If $T$ is a $1$-tree then

\begin{equation}\label{T1}
\left|\sum_{Q\in T}
\langle f_2, \phi_{Q_1}\rangle
\left< \phi_{Q_1} \chi_{\{x:N(x)\in \omega_{Q_2}\}}
,\sum_{P\in\P; \omega_{Q_1}\subseteq\omega_{P_2}}
\langle f_1, \phi_{P_1}\rangle\phi_{P_1}
\chi_{\{x:N(x)\in \omega_{P_2}\}}f_3\right>\right|\lesssim
\end{equation}

$$\size_1((\langle f_2,\phi_{Q_1}\rangle)_{Q\in T})
\size'_{T,\P}(f_1,f_3)|I_T|.$$

If $T$ is a $2$-tree then

\begin{equation}\label{T2'}
\left|\sum_{Q\in T}
\langle f_2, \phi_{Q_1}\rangle
\left<\phi_{Q_1} \chi_{\{x:N(x)\in \omega_{Q_2}\}}
,\sum_{P\in\tilde{T}; \omega_{Q_1}\subseteq\omega_{P_2}}
\langle f_1, \phi_{P_1}\rangle\phi_{P_1}
\chi_{\{x:N(x)\in \omega_{P_2}\}}f_3\right>\right|\lesssim 
\end{equation}

$$
\lesssim\size''_{T,\P}(f_2,f_1)
\size((\langle f_3\chi_{\{x:N(x)\in\omega_{Q_2}\}},
\phi_{Q_1}\rangle )_{Q\in T})|I_T|
$$
and also

\begin{equation}\label{T2''}
\left|\sum_{Q\in T}
\langle f_2, \phi_{Q_1}\rangle
\left< \phi_{Q_1} \chi_{\{x:N(x)\in \omega_{Q_2}\}}
,\sum_{P\in\tilde{T}^c; \omega_{Q_1}\subseteq\omega_{P_2}}
\langle f_1, \phi_{P_1}\rangle\phi_{P_1}
\chi_{\{x:N(x)\in \omega_{P_2}\}}f_3\right>\right|\lesssim 
\end{equation}

$$ \size_1((\langle f_2,\phi_{Q_1}\rangle)_{Q\in T})
\size'_{T,\P}(f_1,f_3)|I_T|
$$
where $\tilde{T}^c$ is the set of all quartiles $P\in\P$ so that
$P_2\leq Q_{T,2}$ but $P$ does not belong to $\tilde{T}$ ($\tilde{T}$
was defined in Definition \ref{size''energy''-def}).

\end{lemma}

\begin{proof}
Let $\cal{J}$ be the collection of all maximal dyadic intervals $J$ so that
$3J$ does not contain any $I_Q$ for $Q\in T$. We observe that $\cal{J}$
is a partition of the real line $\R$.

\underline{Case $1$: $T$ is a $1$-tree}

In this situation, one can estimate the left hand side of (\ref{T1}) by

$$\left\|
\sum_{Q\in T}
\langle f_2, \phi_{Q_1}\rangle
\phi_{Q_1} \chi_{\{x:N(x)\in \omega_{Q_2}\}}
C^c_{walsh,Q,\P}(f_1)
f_3\right\|_1\leq$$

$$\sum_{J\in\cal{J}}\left\|
\sum_{Q\in T}
\langle f_2, \phi_{Q_1}\rangle
\phi_{Q_1} \chi_{\{x:N(x)\in \omega_{Q_2}\}}
C^c_{walsh,Q,\P}(f_1)
f_3\right\|_{L^1(J)}.$$
By using the fact that $T$ is a $1$-tree and from Lemma \ref{biest-trick}
one can see that the last sum is actually equal to

$$\sum_{J\in\cal{J}}\left\|
\sum_{Q\in T}
\langle f_2, \phi_{Q_1}\rangle
\phi_{Q_1} \chi_{\{x:N(x)\in \omega_{Q_2}\}}
C_{walsh,\P}(f_1)
f_3\right\|_{L^1(J)}=$$

$$\sum_{J\in\cal{J}}\sum_{Q\in T}
\int_{\{x\in J\cap I_Q/N(x)\in\omega_{Q_2}\}}
|\langle f_2, \phi_{Q_1}\rangle
\phi_{Q_1}||C_{walsh,\P}(f_1)f_3|\,dx.$$
But this can be majorized by

\begin{equation}\label{T1*}
\size_1((\langle f_2,\phi_{Q_1}\rangle)_{Q\in T})
\sum_{J\in\cal{J}}\sum_{Q\in T}
\int_{\{x\in J\cap I_Q/N(x)\in\omega_{Q_2}\}}
|C_{walsh,\P}(f_1)f_3|\,dx.
\end{equation}
Let us also observe that the $J$ intervals which are relevant for our
summation, are those for which $J\subseteq 3 I_T$, otherwise the
corresponding terms in (\ref{T1*}) are zero.

In order to continue our estimates, we need to make certain
geometric observations about the sets 
$\{x\in J\cap I_Q/N(x)\in\omega_{Q_2}\}$. 

Fix $J\in\cal{J}$
so that $J\subseteq 3 I_T$ and pick a dyadic interval $J'$ so that
$|J'|=2|J|$ and $J\subseteq J'$.
By the maximality of $J$, $3J'$ contains an interval $I_Q$
for some $Q\in T$. We then chose $Q'(=Q'(J))$ with $|I_{Q'}|=|J'|$
and $Q_1\leq Q'_1\leq Q_T$ and observe that for any $Q\in T$
$\{x\in J\cap I_Q/N(x)\in\omega_{Q_2}\}\subseteq
\{x\in I_{Q'(J)}/N(x)\in\omega_{Q'(J)_2}\}$. Moreover, if $Q$ and $Q'$ are
in $T$ and have different scales then the sets
$\{x\in J\cap I_Q/N(x)\in\omega_{Q_2}\}$ and
$\{x\in J\cap I_{Q'}/N(x)\in\omega_{Q'_2}\}$ are disjoint. In particular,
this implies that (\ref{T1*}) can be majorized by

$$\size_1((\langle f_2,\phi_{Q_1}\rangle)_{Q\in T})
\sum_{J\in\cal{J}}
\int_{\{x\in I_{Q'(J)}/N(x)\in\omega_{Q'(J)_2}\} }
|C_{walsh,\P}(f_1)f_3|\,dx=$$

$$\size_1((\langle f_2,\phi_{Q_1}\rangle)_{Q\in T})
\sum_{J\in\cal{J}}
\int_{I_{Q'(J)} }
|C^c_{walsh,Q'(J),\P}(f_1)f_3
\chi_{\{x\in I_{Q'(J)}/N(x)\in\omega_{Q'(J)_2}\} }|\,dx$$
by using again Lemma \ref{biest-trick}. In the end we just observe that
this is smaller than

$$\size_1((\langle f_2,\phi_{Q_1}\rangle)_{Q\in T})
\size'_{T,\P}(f_1,f_3)\sum_{J}|I_{Q'(J)}|\lesssim$$

$$\size_1((\langle f_2,\phi_{Q_1}\rangle)_{Q\in T})
\size'_{T,\P}(f_1,f_3)\sum_{J}|J|\lesssim 
\size_1((\langle f_2,\phi_{Q_1}\rangle)_{Q\in T})
\size'_{T,\P}(f_1,f_3)|I_T|$$
which is exactly (\ref{T1}).

\underline{Case $2$: $T$ is a $2$-tree}

As before, we can majorize the sum of the left hand sides of 
(\ref{T2'}) and (\ref{T2''}) by

$$\sum_J\left\|
\sum_{Q\in T}
\langle f_2, \phi_{Q_1}\rangle
\phi_{Q_1} \chi_{\{x:N(x)\in \omega_{Q_2}\}}
C^c_{walsh,Q,\tilde{T}}(f_1)
f_3\right\|_{L^1(J)}+$$

\begin{equation}\label{T2*}
\sum_J\left\|
\sum_{Q\in T}
\langle f_2, \phi_{Q_1}\rangle
\phi_{Q_1} \chi_{\{x:N(x)\in \omega_{Q_2}\}}
C^c_{walsh,Q,\tilde{T}^c}(f_1)
f_3\right\|_{L^1(J)}= I+II.
\end{equation}
To estimate the first term, we observe that since all the $P's$ are
in $\tilde{T}$, one can majorize the function

$$\sum_{Q\in T}
\langle f_2, \phi_{Q_1}\rangle
\phi_{Q_1} \chi_{\{x:N(x)\in \omega_{Q_2}\}}
C^c_{walsh,Q,\tilde{T}}(f_1)$$
when restricted to a certain $J\in\cal{J}$ pointwise by

$$C\sup_{J\subseteq I}
\frac{1}{|I|}\left|\int_I
\sum_{Q\in T}
\langle f_2, \phi_{Q_1}\rangle
\phi_{Q_1} 
\sum_{P\in\tilde{T}}
\langle f_1, \phi_{P_1}\rangle
\phi_{P_1}\,dx\right|.$$
Using this estimate and the geometric observations discussed earlier,
one can majorize term $I$ in (\ref{T2*}) by

$$\sum_J
\size((\langle f_3\chi_{\{x:N(x)\in\omega_{Q_2}\}},
\phi_{Q_1}\rangle)_{Q\in T})|J|
\sup_{J\subseteq I}
\frac{1}{|I|}\left|\int_I
\sum_{Q\in T}
\langle f_2, \phi_{Q_1}\rangle
\phi_{Q_1} 
\sum_{P\in\tilde{T}}
\langle f_1, \phi_{P_1}\rangle
\phi_{P_1}\,dx\right|\lesssim$$

$$\size((\langle f_3\chi_{\{x:N(x)\in\omega_{Q_2}\}},
\phi_{Q_1}\rangle)_{Q\in T}) 
\left\|\M (\Pi_T(\sum_{Q\in T}
\langle f_2, \phi_{Q_1}\rangle
\phi_{Q_1},
\sum_{P\in\tilde{T}}
\langle f_1, \phi_{P_1}\rangle
\phi_{P_1}))\right\|_1\lesssim$$

$$\size((\langle f_3\chi_{\{x:N(x)\in\omega_{Q_2}\}},
\phi_{Q_1}\rangle)_{Q\in T})
\left\|\Pi_T(\sum_{Q\in T}
\langle f_2, \phi_{Q_1}\rangle
\phi_{Q_1},
\sum_{P\in\tilde{T}}
\langle f_1, \phi_{P_1}\rangle
\phi_{P_1})\right\|_{H^1}\lesssim$$

$$
\size''_{T,\P}(f_2,f_1)
\size((\langle f_3\chi_{\{x:N(x)\in\omega_{Q_2}\}},
\phi_{Q_1}\rangle )_{Q\in T})|I_T|,$$
where $\M (f)$ is the maximal operator defined by

$$\M (f)(x):=\sup_{x\in I}\frac{1}{|I|}\left|\int_I f(y)\,dy\right|.$$
(We used in the above the maximal function characterization of the Hardy
space $H^1$ (see \cite{stein}) and also the fact that $T$ is a $2$-tree).

To estimate the term $II$ in (\ref{T2*}), we first observe by using
Lemma \ref{biest-trick} that

$$\sum_{Q\in T}
\langle f_2, \phi_{Q_1}\rangle
\phi_{Q_1} \chi_{\{x:N(x)\in \omega_{Q_2}\}}
C^c_{walsh,Q,\tilde{T}^c}(f_1)=
\sum_{Q\in T}
\langle f_2, \phi_{Q_1}\rangle
\phi_{Q_1} \chi_{\{x:N(x)\in \omega_{Q_2}\}}
C_{walsh,\tilde{T}^c}(f_1).$$
Then we also remark that the function

$$\sum_{Q\in T}
\langle f_2, \phi_{Q_1}\rangle
\phi_{Q_1} \chi_{\{x:N(x)\in \omega_{Q_2}\}}$$
when restricted to an interval $J\in\cal{J}$ is pointwise smaller
than 

$$\sup_{J\subseteq I}
\frac{1}{|I|}\left|\int_I
\sum_{Q\in T}
\langle f_2, \phi_{Q_1}\rangle
\phi_{Q_1}\,dx\right|.$$
Using these two facts one can argue as before and estimate the term $II$ by

$$\size'_{T,\P}(f_1,f_3)
\left\|\M (\sum_{Q\in T}
\langle f_2, \phi_{Q_1}\rangle
\phi_{Q_1})\right\|_{L^1(3I_T)}\lesssim$$

$$\size'_{T,\P}(f_1,f_3)|I_T|^{1/2}
\left\|\sum_{Q\in T}
\langle f_2, \phi_{Q_1}\rangle
\phi_{Q_1})\right\|_2\lesssim
\size_1((\langle f_2,\phi_{Q_1}\rangle)_{Q\in T})
\size'_{T,\P}(f_1,f_3)|I_T|.$$
The proof is now complete.

\end{proof}

To extend this summation over $T$ to a summation over $\P$ we would like to partition $\P$ into trees $T$ for which one has control over $\sum_T |I_T|$.  
This will be accomplished by several decomposition lemmas. The first
one is well known (see \cite{mtt:walshbiest} for instance, for a proof).

\begin{proposition}\label{decomp}  Let $j = 1,2,3$, 
$\P'$ be a subset of $\P$, $n \in \Z$, $f$ be a function and suppose that
$$ \size_j((\langle f,\phi_{P_j}\rangle)_{P \in \P'} ) \leq 2^{-n} 
\energy_j((\langle f,\phi_{P_j}\rangle)_{P \in \P} ).$$
Then we may decompose $\P' = \P'' \cup \P'''$ such that
\be{size-lower}
\size_j((\langle f,\phi_{P_j}\rangle)_{P \in \P''} ) \leq 2^{-n-1} 
\energy_j((\langle f,\phi_{P_j}\rangle)_{P \in \P} )
\end{equation}
and that $\P'''$ can be written as the disjoint union of trees $\T$ such that
\be{tree-est}
\sum_{T \in \T} |I_T| \lesssim 2^{2n}.
\end{equation}
\end{proposition}
By iterating this proposition one obtains 
(see again \cite{mtt:walshbiest})

\begin{corollary}\label{decomp-cor}  There exists a partition
$$ \P = \bigcup_{n \in \Z} \P_n$$
where for each $n \in \Z$ and $j = 1,2,3$ we have
$$ \size_j((\langle f,\phi_{P_j}\rangle)_{P \in \P_n} ) 
\leq 
\min\left(2^{-n} \energy_j((\langle f,\phi_{P_j}\rangle)_{P \in \P} ) , 
\size_j((\langle f,\phi_{P_j}\rangle)_{P \in \P} )\right).$$
Also, we may cover $\P_n$ by a collection $\T_n$ of trees such that
$$ \sum_{T \in \T_n} |I_T| \lesssim 2^{2n}.$$
\end{corollary}
The next lemma together with its corrolary are also known (see \cite{laceyt3}).

\begin{proposition}\label{decompc}
 Let  
$\Q'$ be a subset of $\Q$, $n \in \Z$, $f$ be a function and suppose that
$$ \size((\langle f\chi_{\{x:N(x)\in\omega_{Q_2}\}},\phi_{Q_1}
\rangle)_{Q\in\Q'}) \leq 2^{-n} 
\energy((\langle f\chi_{\{x:N(x)\in\omega_{Q_2}\}},\phi_{Q_1}
\rangle)_{Q\in\Q})
.$$
Then we may decompose $\Q' = \Q'' \cup \Q'''$ such that
\be{size-lowerc}
\size((\langle f\chi_{\{x:N(x)\in\omega_{Q_2}\}},\phi_{Q_1}
\rangle)_{Q\in\Q''})
\leq 2^{-n-1}
\energy((\langle f\chi_{\{x:N(x)\in\omega_{Q_2}\}},\phi_{Q_1}
\rangle)_{Q\in\Q})
\end{equation}
and that $\Q'''$ can be written as the disjoint union of trees $\T$ such that
\be{tree-estc}
\sum_{T \in \T} |I_T| \lesssim 2^{n}.
\end{equation}
\end{proposition}

\begin{corollary}\label{decomp-corc}  There exists a partition
$$ \Q= \bigcup_{n \in \Z} \Q_n$$
where for each $n \in \Z$  we have
$$ \size((\langle f\chi_{\{x:N(x)\in\omega_{Q_2}\}},\phi_{Q_1}
\rangle)_{Q\in\Q_n})
\leq $$
$$
\min\left(2^{-n} 
\energy((\langle f\chi_{\{x:N(x)\in\omega_{Q_2}\}},\phi_{Q_1}
\rangle)_{Q\in\Q}),
\size((\langle f\chi_{\{x:N(x)\in\omega_{Q_2}\}},\phi_{Q_1}
\rangle)_{Q\in\Q}) \right)
.$$
Also, we may cover $\Q_n$ by a collection $\T_n$ of trees such that
$$ \sum_{T \in \T_n} |I_T| \lesssim 2^{n}.$$
\end{corollary}

We also need

\begin{proposition}\label{decomp'}
 Let  
$\Q'$ be a subset of $\Q$, $n \in \Z$, $f,g$ be two functions and suppose that
$$ \size'_{\Q',\P}(f,g) \leq 2^{-n} 
\energy'_{\Q,\P}(f,g).$$
Then we may decompose $\Q' = \Q'' \cup \Q'''$ such that
\be{size-lower'}
\size'_{\Q'',\P}(f,g) \leq 2^{-n-1}\energy'_{\Q,\P}(f,g) 
\end{equation}
and that $\Q'''$ can be written as the disjoint union of trees $\T$ such that
\be{tree-est'}
\sum_{T \in \T} |I_T| \lesssim 2^{n}.
\end{equation}
\end{proposition}

\begin{proof}
First, let us denote by $\Q'_{heavy}$ the set of all quartiles $Q\in\Q'$
so that

$$\size'_{\{Q\},\P}(f,g)>2^{-n-1}\energy'_{\Q,\P}(f,g).$$
Clearly, by Definition \ref{size'energy'-def}, for every such a quartile
$Q$ there exists a quartile $Q'(Q)\in\Q$ with $Q_{12}<Q'(Q)_{12}$ and so that

$$
\frac{1}{|I_{Q'(Q)}|}
\int_{I_{Q'(Q)}}
|C^c_{walsh,Q'(Q),\P}(F)G|
\chi_{\{x:N(x)\in \omega_{Q'(Q)_1}\cup\omega_{Q'(Q)_2}\}}\,dx
>2^{-n-1}\energy'_{\Q,\P}(f,g).$$
Let then denote by $\D$ the set of all such quartiles $Q'(Q)$
so that $Q'(Q)_{12}$ are maximal with respect to the partial order ``$<$''
previously defined. Because of this maximality, the set $\D$
is a set of quartiles $Q$ so that the corresponding set of sub-bitiles
$Q_{12}$ is a set of disjoint bitiles. Then we define the set
$\T$ to be the collection of all trees in $\Q'$ with top in
$\D$, $\Q'''$ to be union of all tiles in $\T$ and $\Q'':=\Q'\setminus\Q'''$. 
Then (\ref{size-lower'}) follows
by construction while (\ref{tree-est'}) follows from the inequalities

$$\sum_{T\in\T}|I_T|\lesssim
2^n\energy'_{\Q,\P}(f,g)^{-1}\sum_{Q'(Q)\in\D}
\int_{I_{Q'(Q)}}
|C^c_{walsh,Q'(Q),\P}(F)G|
\chi_{\{x:N(x)\in \omega_{Q'(Q)_1}\cup\omega_{Q'(Q)_2}\}}\,dx$$
$$\lesssim 2^n.$$
This completes the proof.
\end{proof}

By iterating this proposition one obtains

\begin{corollary}\label{decomp-cor'}  There exists a partition
$$ \Q= \bigcup_{n \in \Z} \Q_n$$
where for each $n \in \Z$  we have
$$ \size'_{\Q_n,\P}(f,g) 
\leq 
\min(2^{-n} \energy'_{\Q,\P}(f,g) , \size'_{\Q,\P}(f,g))
.$$
Also, we may cover $\Q_n$ by a collection $\T_n$ of trees such that
$$ \sum_{T \in \T_n} |I_T| \lesssim 2^{n}.$$
\end{corollary}
Finally, we need
\begin{proposition}\label{decomp''}
 Let  
$\Q'$ be a subset of $\Q$, $n \in \Z$, $f,g$ be two functions and suppose that
$$ \size''_{\Q',\P}(f,g) \leq 2^{-n} 
\energy''_{\Q,\P}(f,g).$$
Then we may decompose $\Q' = \Q'' \cup \Q'''$ such that
\be{size-lower''}
\size''_{\Q'',\P}(f,g) \leq 2^{-n-1}\energy''_{\Q,\P}(f,g) 
\end{equation}
and that $\Q'''$ can be written as the disjoint union of trees $\T$ such that
\be{tree-est''}
\sum_{T \in \T} |I_T| \lesssim 2^{2n}.
\end{equation}
\end{proposition}

\begin{proof}
The idea is to initialize $\Q''$ to equal $\Q'$, and remove trees from $\Q''$ one by one (placing them into $\Q'''$) until \eqref{size-lower''} is satisfied.

Since $\Q$ has finite cardinality, this procedure will terminate after a finite
number of steps.

We assume by pigeonholing that we only have quartiles $Q$ such that 
the lengths of $I_Q$ are all even (or all odd) powers of $2$.

We describe the tree selection algorithm.  
We shall need four collections $\T'_{i}, \T''_{i}$ of trees, 
where $i \neq 1$; we initialize all four collections to be empty.  

Suppose that we can find an $i \neq 1$ and a quartile $Q^0 \in \Q''$ such that
if we denote by $T= \{ Q\in \Q'': Q_i < Q^0_i \}$, one has

\be{plus}
\|
\sum_{Q\in T}
\langle f,\phi_{Q_1}\rangle\phi_{Q_1}
\sum_{P\in\tilde{T};\omega_{Q_1}\subseteq\omega_{P_2}}
\langle g,\phi_{P_1}\rangle\phi_{P_1}
\|_1\geq 2^{-n-3/2}\energy''_{\Q,\P}(f,g)|I_{Q^0}|.
\end{equation}
We may assume that $Q^0_i$ is maximal with respect to this property and the tile partial order $<$.  Having assumed this maximality, we may then assume that 
$\xi_{Q^0}$ is minimal, where $\xi_{Q^0}$ is the center of $\omega_{Q^0}$.

We then place the $i$-tree
$$ \{ Q \in \Q'': Q_i < Q^0_i \}$$
with top $Q^0$ into the collection $\T'_i$, and then remove all the quartiles in this tree from $\Q''$.  We then place the $1$-tree
$$ \{ Q \in \Q'': Q_1 < Q^0_1 \}$$
with top $Q^0$ into the collection $\T''_i$, and then remove all the quartiles in this tree from $\Q''$.

We then repeat this procedure until there are no further 
quartiles $Q^0 \in \Q''$ which obey \eqref{plus}.

After completing this algorithm, none of the tiles $Q^0$ in $\Q''$ 
will obey  \eqref{plus}, so that \eqref{size-lower''} holds for 
all $i$-trees in $\Q''$.  (If the tree does not contain its top, 
we can break it up as the disjoint union of trees which do).  We then set 
$\T := \bigcup_{i \neq 1} \T'_{i} \cup \T''_{i}$
and $\Q' := \bigcup_{T \in \T} T$.  

It remains to prove \eqref{tree-est''}.  Since the trees in $\T''_{i}$ have the same tops as those in $\T'_{i}$ it suffices to prove the estimate for 
$\T'_i$.  
Fix $i>1$.  The key geometric observation is that the tiles
$$ \{ Q_1:  Q \in T \hbox{ for some } T \in \T'_i \}$$
are all pairwise disjoint.  Indeed, suppose that there existed 
$Q \in T \in T'_i$ and $Q' \in T' \in T'_i$ such that $Q_1 \neq Q'_1$ and 
$Q_1 \cap Q'_1
\neq \emptyset$.  Without loss of generality we may assume that 
\be{G1}
I_Q \supsetneq I_{Q'}
\end{equation}
so that
$$\omega_{Q_1} \subsetneq \omega_{Q'_1}.$$  
From the nesting of dyadic intervals, and from the assumption that
two different scales differ at leats by a factor of $4$, this implies that 
$$\omega_{Q_i} \subsetneq \omega_{Q'_1}.$$
Since $\T'_i$ consists entirely of $i$-trees, we have 
$\omega_{Q_{T,i}} \subset \omega_{Q_i}$, thus
\be{G2}
\omega_{Q_{T,i}} \subsetneq \omega_{Q'_1}.
\end{equation}
On the other hand, since $T'$ is an $i$-tree, we have
$$\omega_{Q_{T',i}} \subseteq \omega_{Q'_i}.$$
Using our selection algorithm,
we thus see that $\omega_{Q_{T,i}}$ and $\omega_{Q_{T',i}}$ 
are disjoint and that
$$ \xi_{Q_{T,i}} < \xi_{Q_{T',i}}.$$
Since we chose our trees $T$ in $\T'_i$ so that $\xi_{Q_{T,i}}$ was minimized,
 this implies that $T$ was selected earlier than $T'$.  On the other hand,
from \eqref{G2} and the nesting of dyadic intervals we have
$$\omega_{Q_{T,1}} \subsetneq \omega_{Q'_1}$$
which implies from \eqref{G1} that
$$Q'_1 < Q_{T,1}.$$
Thus $Q'$ would have been selected for a tree in $\T''_i$ at the same time that $T$ was selected for $\T'_i$.   But this contradicts the fact that $Q'$ is 
part of $T'$, and therefore selected at a later time for $\T'_i$.  
This establishes the pairwise disjointness of the $Q_1$.

Similarly, one can also prove that the tiles
$$ \{ P_1:  P \in \tilde{T} \hbox{ for some } T \in \T'_i \}$$
are all pairwise disjoint.

By using (\ref{plus}), we deduce that for any $T\in\T'_i$ one has

$$|I_T|\lesssim
2^{2n}(\energy''_{\Q,\P}(f,g))^{-2}
\|
\sum_{Q\in T}
\langle f,\phi_{Q_1}\rangle\phi_{Q_1}
\sum_{P\in\tilde{T};\omega_{Q_1}\subseteq\omega_{P_2}}
\langle g,\phi_{P_1}\rangle\phi_{P_1}
\|_2^2.$$
From this and the disjointness of $T$'s and $\tilde{T}$'s, it follows that

$$ \sum_{T \in \T'_i} |I_T| \lesssim 2^{2n}$$
as wanted.

\end{proof}

Once again, by iterating the above proposition we obtain

\begin{corollary}\label{decomp-cor''}  There exists a partition
$$ \Q= \bigcup_{n \in \Z} \Q_n$$
where for each $n \in \Z$  we have
$$ \size''_{\Q_n,\P}(f,g) 
\leq 
\min(2^{-n} \energy''_{\Q,\P}(f,g) , \size''_{\Q,\P}(f,g))
.$$
Also, we may cover $\Q_n$ by a collection $\T_n$ of trees such that
$$ \sum_{T \in \T_n} |I_T| \lesssim 2^{2n}.$$
\end{corollary}

\section{Proof of Propositions \ref{babstract} and \ref{moreabstract}}

It remains to prove Propositions \ref{babstract} and \ref{moreabstract}.
Since Proposition \ref{babstract} is easier, we shall only present
the proof of Proposition \ref{moreabstract}. In the end,
 we will briefly explain
how one can also prove Proposition \ref{babstract} by using some of the same 
ideas.

Fix $\P,\Q$, $f_1,f_2,f_3$ and $\alpha_j ,\beta_j ,\gamma_j\in (0,1)$
as in the hypothesis of Proposition \ref{moreabstract}. By applying
Lemma \ref{decomp} to $f_2$, one obtains a partition

$$\Q=\bigcup_{k\in\Z}\Q^2_k$$
satisfying the conditions in that lemma. In particular, one can decompose
each $\Q^2_k$ as a union of trees in $\T^2_k$. Then, by applying
Lemma \ref{decompc} to $f_3$, one obtains again a partition

$$\Q=\bigcup_{l\in\Z}\Q^3_l$$
satisfying the conditions of that lemma. In particular, one can decompose
each $\Q^3_l$ as a union of trees in $\T^3_l$. Similarly,
by applying Lemma \ref{decomp'} to the pair $(f_1,f_3)$, one obtains
a partition

$$\Q=\bigcup_{m\in\Z}\Q^{13}_m$$
satisfying the conditions of that lemma and in particular one can decompose
each $\Q^{13}_m$ as a union of trees in $\T^{13}_m$.

Finally, by applying this time Lemma \ref{decomp''} to the pair
$(f_2,f_1)$ one obtains another decomposition

$$\Q=\bigcup_{n\in\Z}\Q^{21}_n$$
satisfying the conditions of that lemma. In particular, one can
decompose as before each $\Q^{21}_n$ as a union of trees in $\T^{21}_n$.

Using the above decompositions, one can write our form
$\Lambda''(f_1,f_2,f_3)$ as

$$\Lambda''(f_1,f_2,f_3)=
\sum_{Q\in \Q}
\langle f_2, \phi_{Q_1}\rangle
\phi_{Q_1} \chi_{\{x:N(x)\in \omega_{Q_2}\}}
C^c_{walsh,Q,\P}(f_1)
f_3=$$

$$\sum_{k,l,m,n}
\sum_{Q\in \Q^2_k\cap\Q^3_l\cap\Q^{13}_m\cap\Q^{21}_n}
\langle f_2, \phi_{Q_1}\rangle
\phi_{Q_1} \chi_{\{x:N(x)\in \omega_{Q_2}\}}
C^c_{walsh,Q,\P}(f_1)
f_3=$$

\begin{equation}\label{a}
\sum_{k,l,m,n}
\sum_{T\in\T^{k,l,m,n}}
\sum_{Q\in T}
\langle f_2, \phi_{Q_1}\rangle
\phi_{Q_1} \chi_{\{x:N(x)\in \omega_{Q_2}\}}
C^c_{walsh,Q,\P}(f_1)
f_3
\end{equation}
where $\T^{k,l,m,n}$ is a collection of $i$-trees of tiles in $\Q$
with the property that every tree $T\in\T^{k,l,m,n}$ is actually
a subtree in $\T^2_k$, $\T^3_l$, $\T^{13}_m$ and $\T^{12}_n$. 
We can also naturally decompose 
$\T^{k,l,m,n}=\T_1^{k,l,m,n}\cup\T_2^{k,l,m,n}$ where 
$\T_i^{k,l,m,n}$ contains $i$-trees only. Consequently, our sum
(\ref{a}) splits into

$$
\sum_{k,l,m,n}
\sum_{T\in\T_1^{k,l,m,n}}
\sum_{Q\in T}
\langle f_2, \phi_{Q_1}\rangle
\phi_{Q_1} \chi_{\{x:N(x)\in \omega_{Q_2}\}}
C^c_{walsh,Q,\P}(f_1)
f_3+$$

$$\sum_{k,l,m,n}
\sum_{T\in\T_2^{k,l,m,n}}
\sum_{Q\in T}
\langle f_2, \phi_{Q_1}\rangle
\phi_{Q_1} \chi_{\{x:N(x)\in \omega_{Q_2}\}}
C^c_{walsh,Q,\P}(f_1)
f_3=$$

$$
\sum_{k,l,m,n}
\sum_{T\in\T_1^{k,l,m,n}}
\sum_{Q\in T}
\langle f_2, \phi_{Q_1}\rangle
\phi_{Q_1} \chi_{\{x:N(x)\in \omega_{Q_2}\}}
C^c_{walsh,Q,\P}(f_1)
f_3+$$

$$\sum_{k,l,m,n}
\sum_{T\in\T_2^{k,l,m,n}}
\sum_{Q\in T}
\langle f_2, \phi_{Q_1}\rangle
\phi_{Q_1} \chi_{\{x:N(x)\in \omega_{Q_2}\}}
C^c_{walsh,Q,\tilde{T}^c}(f_1)
f_3+$$

$$\sum_{k,l,m,n}
\sum_{T\in\T_2^{k,l,m,n}}
\sum_{Q\in T}
\langle f_2, \phi_{Q_1}\rangle
\phi_{Q_1} \chi_{\{x:N(x)\in \omega_{Q_2}\}}
C^c_{walsh,Q,\tilde{T}}(f_1)
f_3=$$

$$= I+II+III.$$
By using now the above lemmas, the term $I+II$ can be estimated
(after regrouping the terms) by

\begin{equation}\label{b}
\energy_1(f_2)
\energy'(f_1,f_3)
\sum_{k,m}2^{-k}2^{-m}
\sum_{T\in\T^{k,m}}|I_T|
\end{equation}
where $\T^{k,m}$ is a set of trees which are subtrees in $\T^2_k$ and
$\T^{13}_m$ and the summation in (\ref{b}) runs over the indices
$k,m\in\Z$ so that

$$2^{-k}\lesssim\frac{\size_1(f_2)}{\energy_1(f_2)}$$
and

$$2^{-m}\lesssim\frac{\size'(f_1,f_3)}{\energy'(f_1,f_3)}.$$
Moreover, one can estimate the total number of trees in two different ways,
namely

$$\sum_{T\in\T^{k,m}}|I_T|\lesssim 2^{2k},\,\,2^{m}.$$
By interpolating between these two estimates one obtains

\begin{equation}\label{c}
\sum_{T\in\T^{k,m}}|I_T|\lesssim 2^{2a_1 k}2^{a_2 m}
\end{equation}
for every $a_1\in (0,1/2)$, $a_2\in (0,1)$ so that $a_1+a_2=1$.
By inserting (\ref{c}) into (\ref{b}) one obtains after summing over
$k,m$ the majorant

\begin{equation}\label{d}
\size_1(f_2)^{1-2a_1}
\size'(f_1,f_3)^{1-a_2}
\energy_1(f_2)^{2a_1}
\energy'(f_1,f_3)^{a_2}.
\end{equation}
Now, if we set $1-2a_1:=\alpha_1$ and $1-a_2:=\alpha_2$ this bound (\ref{d})
coincides with the first right hand side term in the inequality
of Proposition \ref{moreabstract}.

It remains to estimate the term $III$. We shall do this in two different ways.

Firstly, by applying again the above lemmas, one can estimate it
(after regrouping the tiles) by

\begin{equation}\label{e}
\energy''_1(f_2,f_1)
\energy(f_3)
\sum_{l,n}2^{-l}2^{-n}
\sum_{T\in\T^{l,n}}|I_T|
\end{equation}
where $\T^{l,n}$ is a set of trees which are subtrees in each
$\T^3_l$ and $\T^{21}_n$ and the summation runs over the indices
$l,n\in\Z$ satisfying

$$2^{-l}\lesssim\frac{\size''_1(f_2,f_1)}{\energy''_1(f_2,f_1)}$$
and

$$2^{-n}\lesssim\frac{\size(f_3)}{\energy(f_3)}.$$
As before, the total number of trees can be estimated in two different
ways, namely

$$\sum_{T\in\T^{l,n}}|I_T|\lesssim 2^{2l},\,\,2^{n}$$
and in particular we also get

\begin{equation}\label{f}
\sum_{T\in\T^{l,n}}|I_T|\lesssim 2^{2b_1 l}2^{b_2 n}
\end{equation}
for every $b_1\in (0,1/2)$, $b_2\in (0,1)$ so that $b_1+b_2=1$. By inserting
this estimate (\ref{f}) into (\ref{e}) one obtains after summing over
$l,n$ the majorant

\begin{equation}\label{g}
\size''_1(f_2,f_1)^{1-2b_1}
\size(f_3)^{1-b_2}
\energy''_1(f_2,f_1)^{2b_1}
\energy(f_3)^{b_2}.
\end{equation}
If we set now $1-2b_1=\beta_1$ and $1-b_2=\beta_2$ this bound (\ref{g})
becomes $E_1(f_1,f_2,f_3,(\beta_j)_j)$ in Proposition \ref{moreabstract}.

It thus suffices to show that the term $III$ can also be estimated by
$E_2(f_1,f_2,f_3,(\gamma_j)_j)$, in order to complete our proof.

At first, after regrouping the tiles we rewrite $III$ as

\begin{equation}\label{h}
III=
\sum_{k,l}
\sum_{T\in\T_2^{k,l}}
\sum_{Q\in T}
\langle f_2, \phi_{Q_1}\rangle
\phi_{Q_1} \chi_{\{x/N(x)\in \omega_{Q_2}\}}
C^c_{walsh,Q,\tilde{T}}(f_1)
f_3
\end{equation}
where $\T_2^{k,l}$ is a set of $2$-trees which are subtrees in
$\T^2_k$ and $\T^3_l$.

Fix now $k,l$. The tiles $P$ in (\ref{h}) run inside the set

$$\bigcup_{T\in\T^{k,l}_2}\tilde{T}:=\tilde{\P}^{k,l}.$$
Also, by construction, all these  $\tilde{T}$ $2$-trees are disjoint and
they can be thought of
as being trees of $P$ tiles with tops in $\Q$. 

We should also point out here the general straightforward geometric fact,
that if two trees $T'$ and $T''$ are maximal with respect to
inclusion and they lie inside
the same tree $T'''$, then $I_{T'}\cap I_{T''}=\emptyset$.

Using the disjointness of our $\tilde{T}$ $2$-trees,
the above geometric observation and the estimate 
$\|f_1\|_{\infty}\lesssim 1$, we can naturally
decompose the set $\tilde{\P}^{k,l}$ in the spirit of the above
lemmas as

$$\tilde{\P}^{k,l}=\bigcup_{n'\in\Z}\tilde{\P}^{k,l}_{n'}$$
where for each $n'\in\Z$ we have

$$\sup_{\T'\subseteq \tilde{\P}^{k,l}_{n'}}
|I_{T'}|^{-1/2}
(\sum_{P\in T'}
\langle f_1,\phi_{P_1}\rangle^2)^{1/2}
\lesssim \min (2^{-n'}\energy_1(f_1), 1).$$
Also, we may cover $\tilde{\P}^{k,l}_{n'}$ by a collection of trees
$\tilde{\T}^{k,l}_{n'}$   so that 

$$\sum_{T'\in\tilde{\T}^{k,l}_{n'}}|I_{T'}|
\lesssim \min (2^{2n'}, \sum_{T\in\T_2^{k,l}}|I_T|).
$$
(The $2$-trees $T'$ in the above decomposition are also trees of $P$ tiles
with tops in $\Q$ !).

Now, by using this new splitting, the other decomposition lemmas,
the tree estimate (\ref{T2''}), the definition of 
$\size''(f_2,f_1)$ and the $L^2\times L^2\rightarrow L^1$ boundedness
of paraproducts, one can estimate the absolute value of $III$ by

\begin{equation}\label{i}
\energy_1(f_2)
\energy(f_3)
\energy_1(f_1)
\sum_{k,l,n'} 2^{-k}2^{-l}2^{-n'}
\sum_{T\in\T_2^{k,l,n'}}|I_T|
\end{equation}
where $\T_2^{k,l,n'}$ is a set of $2$-trees $T$ which are subtrees in
$\T^{k,l}_2$ and so that their 
$\tilde{T}$'s are subtrees in $\tilde{\T}^{k,l}_{n'}$. Also, the parameters
$k,l,n'$ in the above summation satisfy the constraints

$$2^{-k}\lesssim\frac{\size_1(f_2)}{\energy_1(f_2)},$$

$$2^{-l}\lesssim\frac{\size(f_3)}{\energy(f_3)}$$ 
and

$$2^{-n'}\lesssim\frac{1}{\energy_1(f_1)}$$
and $\sum_{T\in\T_2^{k,l,n'}}|I_T|$ can be estimated in three different
ways, namely

\begin{equation}
\sum_{T\in\T_2^{k,l,n'}}|I_T|\lesssim 2^{2k},\,\,2^{l},\,\,2^{2n'}.
\end{equation}
By interpolating these inequalities we get

\begin{equation}\label{j}
\sum_{T\in\T_2^{k,l,n'}}|I_T|\lesssim 2^{2c_1 k}2^{c_2 l}2^{2c_3 n'},
\end{equation}
where $c_1,c_3\in (0,1/2)$, $c_2\in (0,1)$ and $c_1+c_2+c_3=1$.
By using (\ref{j}) into (\ref{i}) we obtain after summing over
$k,l,n'$ the majorant

\begin{equation}\label{k}
\size_1(f_2)^{1-2c_1}
\size(f_3)^{1-c_2}
\energy_1(f_2)^{2c_1}
\energy_1(f_1)^{2c_3}
\energy(f_3)^{c_2}.
\end{equation}
In the end, this bound (\ref{k}) becomes equal to
$E_2(f_1,f_2,f_3,(\gamma_j)_j)$ if we set
$1-2c_1=\gamma_1$, $1-c_2=\gamma_3$ and $1-2c_3=\gamma_2$.

This completes the proof of Proposition \ref{moreabstract}.

To prove Proposition \ref{babstract}, one argues in the same way.
Since in this case there is no double summation, the complicated terms
$E_j(f_1,f_2,f_3)$, $j=1,2$ will simply disapear in our previous proof
and what remains is precisely the inequality stated in
Proposition \ref{babstract}.

\end{document}